\documentclass[a4paper,12pt]{article}
\usepackage{amsmath, amsthm, amssymb}
\usepackage{graphicx,subfigure}
\usepackage[geometry]{ifsym}
\usepackage{anysize}
\usepackage{float}
\usepackage{caption}
\marginsize{2.5cm}{2.5cm}{2cm}{2cm}

\newtheorem{theorem}{Theorem}[section]
\newtheorem{lemma}[theorem]{Lemma}
\newtheorem{proposition}[theorem]{Proposition}
\newtheorem{property}[theorem]{Property}
\newtheorem{corollary}[theorem]{Corollary}
\newtheorem{definition}[theorem]{Definition}

\begin{document}
\newcommand{\Eqref}[1]{(\ref{#1})}

\def\CC{{\rm C\kern-.18cm\vrule width.6pt height 6pt depth-.2pt
\kern.18cm}}

\def\NN{{\mathop{{\rm I}\kern-.2em{\rm N}}\nolimits}}

\def\ZZ{{\mathbb{Z}}}

\def\PP{{\mathop{{\rm I}\kern-.2em{\rm P}}\nolimits}}

\def\RR{{\mathop{{\rm I}\kern-.2em{\rm R}}\nolimits}}

\title{Generalized T-splines and VMCR T-meshes}

\author{{Cesare Bracco}$^a$ and {Durkbin Cho}$^b$\footnote{Corresponding author. E-mail: durkbin@dongguk.edu, Phone: +82 2 2260 3197, Fax: +82 2 2290 1380.}%, \\ and {Tae-wan Kim}$^{b}
\\ \\
\small $^a$Department of Mathematics \lq\lq G. Peano\rq\rq, University of Turin\\ \small V. Carlo Alberto 10, Turin 10123, Italy\\
%\small $^b$Department of Naval Architecture \& Ocean Engineering - Seoul National University\\ \small 1, Gwanak-ro, Gwanak-gu, Seoul 151-744, Korea\\
\small $^b$Department of Mathematics, Dongguk University - Seoul \\ \small Pil-dong 3-ga, Jung-gu, Seoul 100-715, Korea\\%\small \texttt{cesare.bracco@yahoo.it}\\
%\and
}\date{}
\maketitle

%******************************* BODY ******************************

%%%%%%%%%%%%%%%%%%%%%%%%%%%%%%%%%%%%%%%%%%%%%
%%%%%%%%%%%%%%%%%%%%%%%%%%%%%%%%%%%%%%%%%%%%%
 
%%%%%%%%%%%%%%%%%%%%%%%%%%%%%%%%%%%%%%%%%%%%%
%%%%%%%%%%%%%%%%%%%%%%%%%%%%%%%%%%%%%%%%%%%%%

\begin{abstract}
The paper considers the extension of the T-spline approach to the Generalized B-splines (GB-splines), a relevant class of non-polynomial splines. The Generalized T-splines (GT-splines) are based both on the framework of classical polynomial T-splines and on the Trigonometric GT-splines (TGT-splines), a particular case of GT-splines. Our study of GT-splines introduces a class of T-meshes (named VMCR T-meshes) for which both the corresponding GT-splines and the corresponding polynomial T-splines are linearly independent. A practical characterization can be given for a sub-class of VMCR T-meshes, which we refer to as weakly dual-compatible T-meshes, which properly includes the class of dual-compatible (equivalently, analysis-suitable) T-meshes for an arbitrary (polynomial) order.
\end{abstract}


\section{Introduction}

\let\thefootnote\relax\footnotetext{Email addresses: cesare.bracco@unito.it (C. Bracco).}

\medskip
In the last years, the introduction of the so-called T-splines and of the spline spaces defined over T-meshes introduced significant advancements for the use of polynomial spline functions in the CAD and CAGD techniques. The main idea of this approach, in the basic case of surface modelling in $\RR^3$, is to free the control points of the surface from the constraint to lie, topologically, on a rectangular grid whose edges intersect only at \lq\lq cross junctions\rq\rq, and allow instead partial lines of control points, which leads to the possibility to have \lq\lq T-junctions\rq\rq  between the edges of the grid. Such a framework gave some important improvements in CAD and CAGD methods: the possibility to locally refine the surfaces, a considerable reduction of the quantity of control points needed, the ability to easily avoid gaps between surfaces to be joined (see, e.g., \cite{sederberg03} and \cite{sederberg04}), just to name a few. All these advantages became even more important in the applications, such as the isogeometric approach for the analysis problems represented by partial differential equations (see, e.g., \cite{cottrell}, \cite{hughes05} and \cite{bazilevs10}).

%\noindent
The T-spline idea has been applied mainly to polynomial splines, while we know that several types of non-polynomial splines are used for certain applications because of their particular properties. For this reason, recently we proposed a generalization of the T-spline approach to the trigonometric GB-splines (see \cite{bracco12b}), a particularly relevant class of non-polynomial splines because of their adaptability and their application to the already mentioned isogeometric analysis (see, e.g., \cite{kvasov99} and \cite{manni11}). Roughly speaking, the GB-splines are a basis of spaces of piecewise functions, locally spanned both by polynomials and by non-polynomial functions, which in the trigonometric case are $\sin(\omega s)$ and $\cos(\omega s)$, with a given frequency $\omega$. Note that these splines can be seen as particular cases of the piecewise Extended Chebyshevian splines (see, e.g., \cite{mazure04}, \cite{mazure11a} and \cite{mazure11b}). GB-splines have been successfully used to construct tensor-product surfaces (see, e.g., \cite{manni11} and references therein) with control points on rectangular grids.

\medskip
%\noindent
In this paper, we will first extend the results in \cite{bracco12b} to any type of GB-splines of {\it arbitrary bi-order $(p,q)$}, so that we can take full advantage of the good features of GB-splines and T-splines. In order to achieve this goal, we will start by presenting the univariate GB-splines and their properties, including a knot insertion formula with necessary conditions, which will be also essential in the study of the linear independence of the GT-spline functions. Then, we will introduce the GT-splines, whose definition (and notations) is based both on the polynomial T-splines (see, e.g., \cite{cho12}, \cite{daveiga} and \cite{li12}) and on the TGT-splines (see \cite{bracco12b}). Similarly to the case of TGT-splines \cite{bracco12b}, we will show that there exists a relation between the GT-splines of bi-order $(p,q)$ and the polynomial T-splines of the same bi-order. The study of their linear independence will lead us to the introduction of the class of VMCR T-meshes (Void Matrix after Column Reduction T-meshes), which guarantee the linear independence of the associated GT-spline and T-spline blending functions of the same bi-order. The basic concept behind VMCR T-meshes involves the idea of column reduction (used in \cite{li12}), and implicitly helped to show in \cite{bracco12b} that {\sl in the case of bi-order $(4,4)$} the well-known analysis-suitable T-meshes are also VMCR T-meshes. In this paper we provide a simple characterization of a sub-class of VMCR T-meshes, which we refer to as {\sl weakly dual-compatible T-meshes}: we will prove that such class strictly includes the one of dual-compatible/analysis-suitable T-meshes (see, e.g., \cite{daveiga} and \cite{li12}) for any bi-order $(p,q)$. Finally, we will present an explicit example of weakly dual-compatible T-meshes which is not dual-compatible/analysis-suitable.

\medskip
%\noindent
The paper is organized as follows. In Section 2 we recall the definition and the basic properties of the univariate GB-splines, and we deal with the conditions needed to get a knot insertion formula. In Section 3, after having recalled the definition of T-mesh, we introduce the GT-splines and we give some properties following directly from their definition. In Section 4 we study the linear independence of the GT-spline blending functions and, more importantly, the classes of VMCR T-meshes and of weakly dual-compatible T-meshes. Finally, Section 5 contains some concluding remarks.

\section{Univariate generalized B-splines}

\subsection{Definition and main properties}

\medskip
Let $n,p\in\NN$, $p\ge 2$, and let ${\bf \Sigma}=\{s_1\le ...\le s_{n+p}\}$ be a non-decreasing knot sequence ({\sl knot vector}); we associate to ${\bf \Sigma}$ two vectors of functions ${\bf \Omega_u}=\{u_1(s),...,u_{n+p-1}(s)\}$ and ${\bf \Omega_v}=\{v_1(s),...,v_{n+p-1}(s)\}$, where, for $i=1,...,n+p-1$, $u_i,v_i$ belong to ${C}^{p-2}[s_i,s_{i+1}]$ and are such that the space $W$ spanned by the derivatives
\begin{equation*}
U_i(s)=\frac{d^{p-2}u_i(s)}{ds^{p-2}}, \quad V_i(s)=\frac{d^{p-2}v_i(s)}{ds^{p-2}}
\end{equation*}
is a Chebyshev space, that is, any function belonging to it has at most one zero in $[s_i,s_{i+1}]$.
Let, for $i=1,...,n+p$, $m_i$ be the multiplicity of $s_i$ in ${\bf \Sigma}$, that is, the cardinality of the set
$$
\{k:\,1\le k\le n+p,\,s_k=s_i \}.
$$
Note that $m_i=m_j$ if $s_i=s_j$. We assume that $1\le m_i\le p$, for $i=1,...,n+p$.
We consider the {\sl generalized spline space} spanned, in each interval $[ s_i, s_{i+1}]$, by $\{ u_i(s), v_i(s), 1, s,..., s^{p-3}\}$ for $p\ge 3$ and by $\{ u_i(s), v_i(s)\}$ for $p=2$. For this space we can define a basis of compactly-supported splines, which are called {\sl Generalized B-splines} (GB-splines). 

\noindent
The definition of such basis is usually given in a recursive fashion, which we briefly recall (see also \cite{kvasov99} and \cite{manni11}). Since we required that the space spanned by $U_i$ and $V_i$, denoted by $W=\langle U_i, V_i\rangle$, is a Chebyshev space, %t contains a unique element which assumes the values $0$ and $1$ at $s_i$ and $s_{i+1}$, respectively. Such an element is positive on $(s_i,s_{i+1})$, since any non-zero element of $\langle U_i, V_i\rangle$ has at most one zero in $[s_i,s_{i+1}]$. Similarly, we can determine another element assuming the values $1$ and $0$ at $s_i$ and $s_{i+1}$, respectively.
it is not restrictive to choose, as generating functions of $W$, two functions $U_i(s)$ and $V_i(s)$ such that
\begin{equation}
U_i(s_i)>0, \quad U_i(s_{i+1})=0, \qquad V_i(s_i)=0, \quad V_i(s_{i+1})>0. \label{an1}
\end{equation}
Since $W$ is a Chebyshev space, $U_i(s)>0$ for any $s\in[s_i,s_{i+1})$ and $V_i(s)>0$ for any $s\in (s_i,s_{i+1}]$. We will call the selected functions $U_i(s)$ and $V_i(s)$ {\sl generating functions associated to $[s_i,s_{i+1}]$}. Then, following \cite{kvasov99}, \cite{manni11}, we can define a basis of compactly-supported spline functions for the generalized spline space in the following way: for $p=2$
\begin{equation}
N_{i}^{(2)}(s)=\begin{cases}
\frac{\displaystyle V_i(s)}{\displaystyle V_i(s_{i+1})},&\text{if $s_i\le s< s_{i+1}$,} \\ 
\frac{\displaystyle U_{i+1}(s)}{\displaystyle U_{i+1}(s_{i+1})},&\text{if $s_{i+1}\le s< s_{i+2}$,} \\ 
0,&\text{otherwise,}
\end{cases}, \label{a1}
\end{equation}
while, for $p\ge 3$,
\begin{equation}
N_{i}^{(p)}(s)=\int_{-\infty}^{s} \big(\delta_{i}^{(p-1)}N_{i}^{(p-1)}(r)-\delta_{i+1}^{(p-1)}N_{i+1}^{(p-1)}(r)\big)dr, \qquad i=1,...,n, \label{a2}
\end{equation}
where
\begin{equation}
\delta_{i}^{(p)}=\Bigg[\int_{-\infty}^{\infty}N_{i}^{(p)}(r)dr\Bigg]^{-1}, \qquad i=1,...,n. \label{a3}
\end{equation}
Moreover, if $N_{i}^{(p)}(s)=0$, we set
\begin{equation*}
\int_{-\infty}^{s}\delta_{i}^{(p)}N_{i}^{(p)}(r) dr=\begin{cases}
1,&\text{$s\ge s_{i+p}$,} \\ 
0,&\text{$s< s_{i+p}$,}
\end{cases}
\end{equation*}

%\noindent
%{\bf Remark.} The assumption that $\langle U_i, V_i\rangle$ is a Chebyshev space is needed for several reasons: it is necessary to prove that the functions defined in \eqref{a1} and \eqref{a2} are a basis of the generalized spline space (see \cite{manni11}) and it guarantees the positivity of such functions. Moreover, it is essential to have a well-posed definition of the basis: however the generating function $U_i(s)$ ($V_i(s)$) belonging to $W$ is chosen, the normalized function $U_i(s)/U_i(s_{i})$ ($V_i(s)/V_i(s_{i+1})$) does not change, since, being $W$ a Chebyshev space, there is a unique element of $W$ taking values $1$ ($0$, resp.) at $s_i$ and $0$ ($1$, resp.) at $s_{i+1}$.

\smallskip
\noindent
The GB-splines have essentially the same properties of the classical polynomial splines.

\smallskip
\begin{property}
The GB-splines satisfy the following properties.
\begin{enumerate}
\item {\bf Continuity:} each $N_{i}^{(p)}$ is $(p-m_j-1)$ times continuously differentiable at the knot $s_j$, where $m_j$ is the multiplicity of $s_j$ in the knot vector $\{s_i,...,s_{i+p}\}$, with $1\le m_i\le p$.
\item {\bf Positivity:} $N_{i}^{(p)}(s)\ge 0$ for $s\in\RR$, $i=1,...,n$ and $p\in\NN$, $p\ge 2$.
\item {\bf Local support:} if $s\notin [s_i,s_{i+p}]$ $N_{i}^{(p)}(s)=0$, $i=1,...,n$ and $p\in\NN$, $p\ge 2$.
\item {\bf Partition of unity:} for $p\ge 3$ and $s\in[s_p,s_{n+1}]$, $\sum_{i=1}^{n} N_{i}^{(p)}(s)=1$.
\item {\bf Linear independence:} for any $p\ge 2$ $N_{1}^{(p)},...,N_{n}^{(p)}$ are linearly independent.
\end{enumerate}
\end{property}

\noindent
Therefore, the GB-splines can be used, just like the polynomial B-splines, to define a GB-spline curve:
\begin{equation}\label{127feb14}
C(t)=\sum_{i=1}^n P_i N_i^{(p)}(s), \qquad s\in[s_1,s_{n+p}],
\end{equation}
where $P_i$ are the control points corresponding to $N_i^{(p)}(s)$. Because of the properties of positivity, local support and partition of unity, the control points in \eqref{127feb14} play the same role of the control points in the polynomial B-spline curves.

\begin{figure}[H] 
\begin{minipage}[H]{.45\textwidth}
\centering
\includegraphics[scale=0.44]{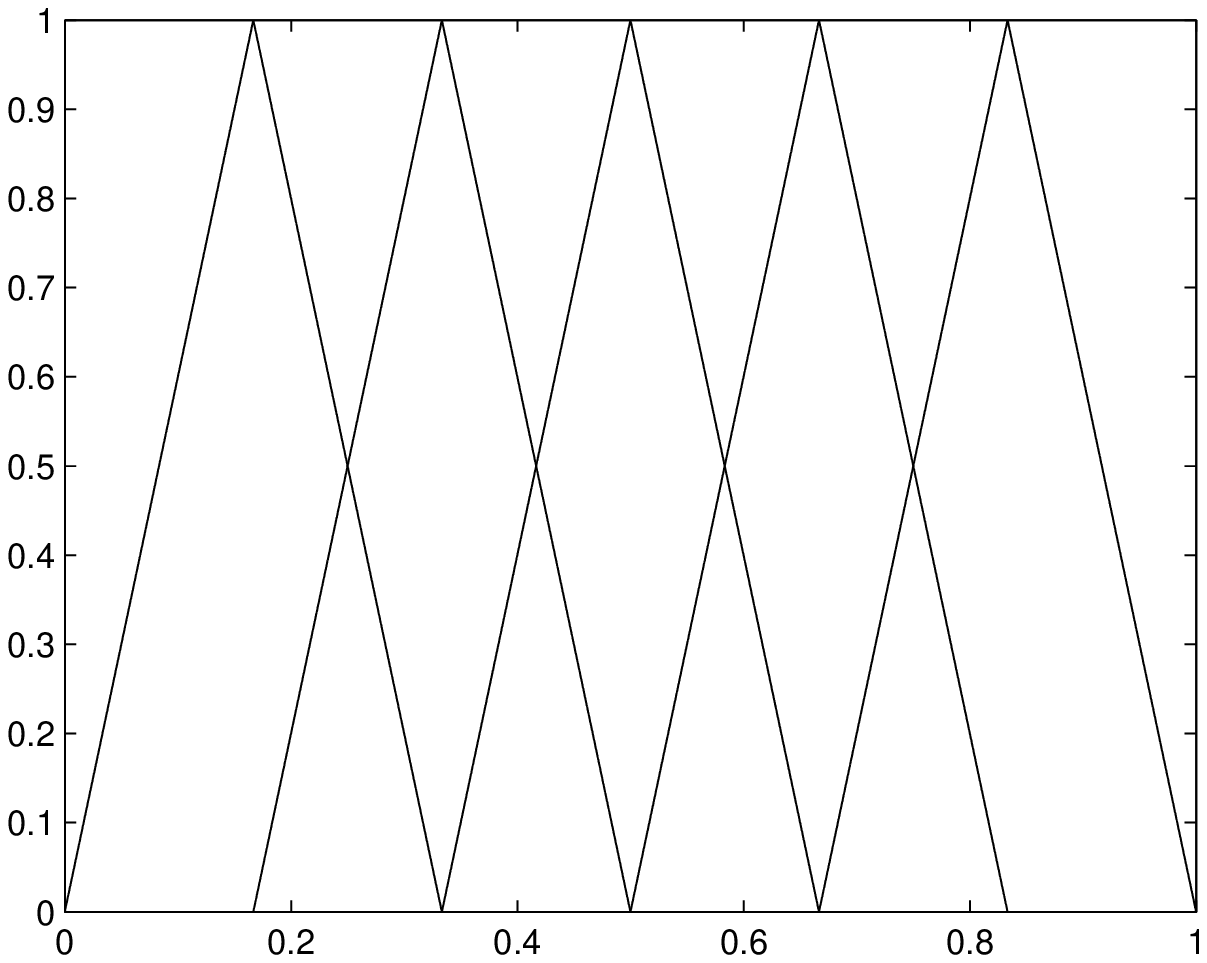}
\caption{\sl (a)}
\end{minipage} 
\hspace{.05\textwidth}
\begin{minipage}[H]{.45\textwidth}
\centering
\includegraphics[scale=0.44]{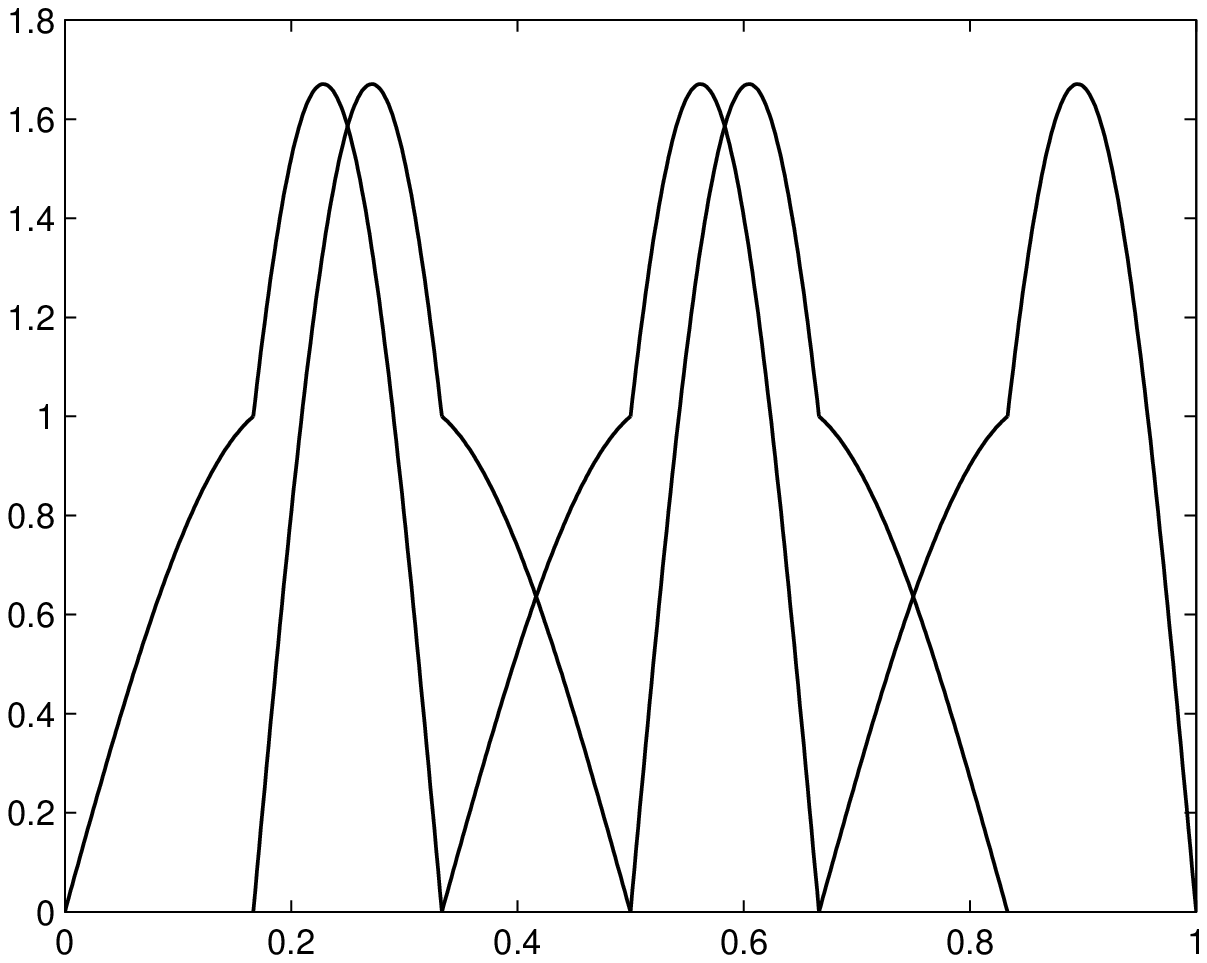}
\caption{\sl (b)}
\end{minipage}
\caption{Figure 1: {\sl A comparison between (a) the polynomial B-splines of order $p=2$ defined on the knot vector ${\bf \Sigma}=\{0,1/6,1/3,1/2,2/3,5/6,1 \}$  and (b) the GB-splines of order $p=2$, locally spanned by the functions $\{u_i(s), v_i(s)\}$, where $u_i(s)=cos(\omega_i s)$ and $v_i(s)=sin(\omega_i s)$ with $\{\omega_i\}_{i=1}^6=\{8,15,8,15,8,15\}$. Note that, since $p=2$, the GB-splines do not form a partition of unity.}}
\end{figure}
\begin{figure}[H] 
\begin{minipage}[H]{.45\textwidth}
\centering
\includegraphics[scale=0.44]{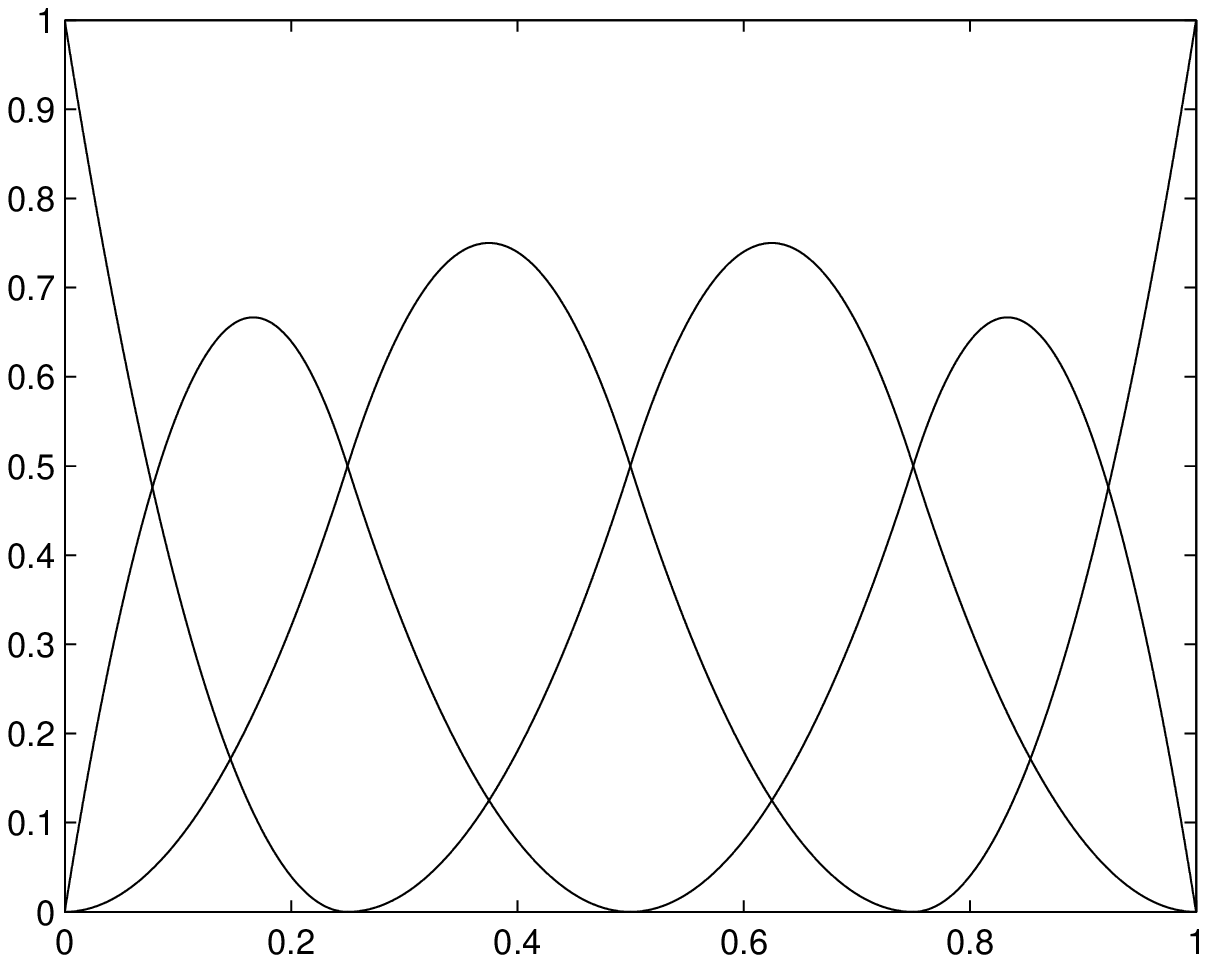}
\caption{\sl (a)}
\end{minipage} 
\hspace{.05\textwidth}
\begin{minipage}[H]{.45\textwidth}
\centering
\includegraphics[scale=0.44]{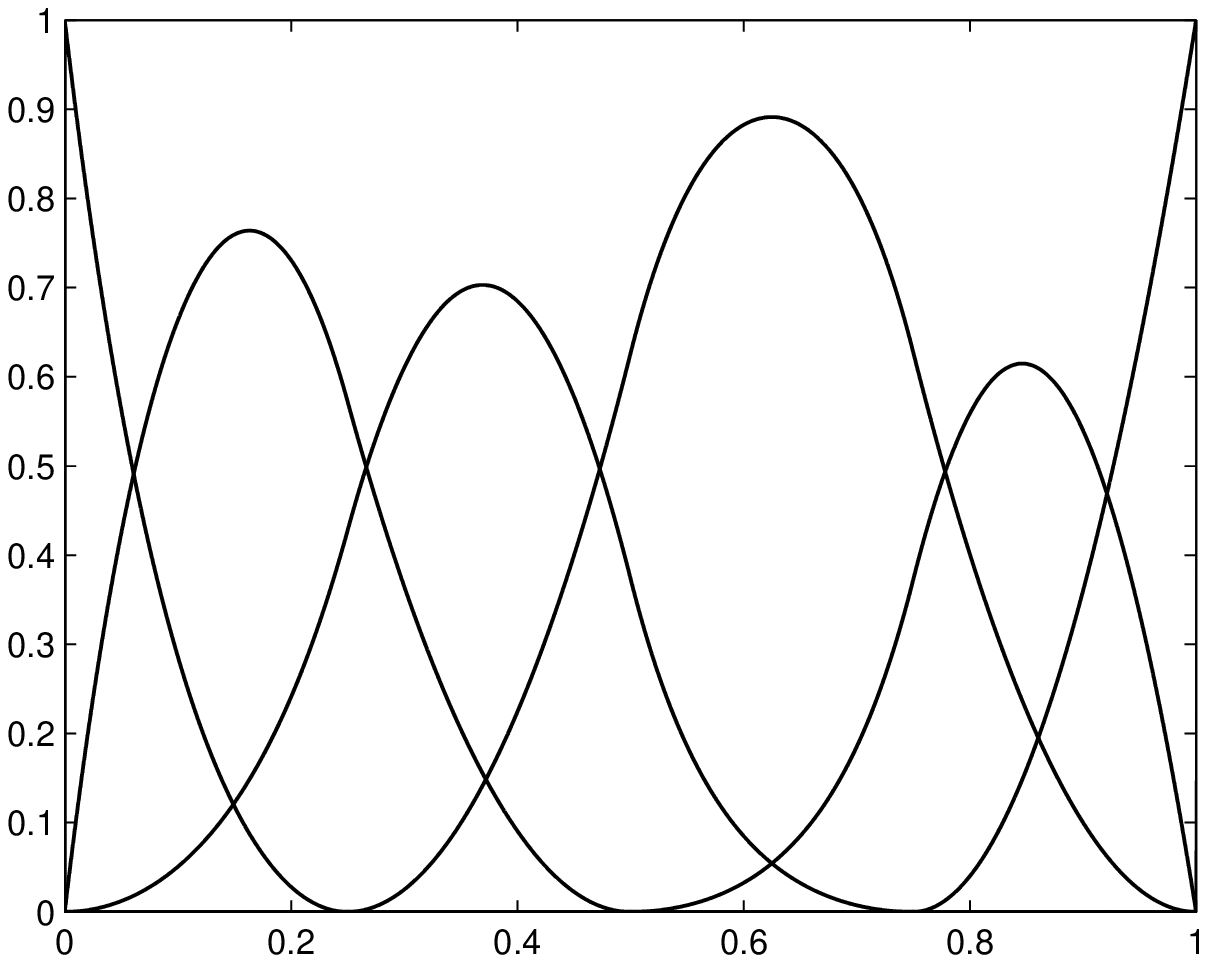}
\caption{\sl (b)}
\end{minipage}
\caption{Figure 2: {\sl A comparison between (a) the polynomial B-splines of order $p=3$ defined on the knot vector ${\bf \Sigma}=\{0,0,0,1/4,1/2,3/4,1,1,1 \}$ and (b) the GB-splines of order $p=3$, locally spanned by the functions $\{1, u_i(s), v_i(s)\}$, where $u_i(s)=cosh(\omega_i s)$ and $v_i(s)=sinh(\omega_i s)$ with $\{\omega_i\}_{i=1}^8=\{1,1,8.5,1.3,12.3,0.5,1,1\}$ for $i=1,...,8$.}}
\end{figure}

\subsection{Knot insertion formula}

\medskip
%\noindent
One of the main reasons to introduce the T-spline approach in the construction of spline surfaces is, as already mentioned, the possibility to apply local refinement techniques. Therefore, it is crucial to have a knot insertion formula in the univariate case. In the case of GB-splines, differently from the polynomial case (see Boehm's seminal work \cite{boehm}), we need to pay attention to the issue that a knot insertion requires two new additional functions and that the space refined by applying the knot insertion formula must contain the original space. The resulting knot insertion rule is stated below. %The formula in \cite{bracco12b} for the trigonometric GB-splines is a particular case of \eqref{a5}; similar, but different, results are available in \cite{wang04} and \cite{wang08} for a particular case of GB-splines and for a slightly different spline space, respectively.

\medskip
\begin{theorem}\label{00a}
Let ${\bf \Sigma}=\{s_1,...,s_{n+p}\}$ be a knot vector, ${\bf \bar \Sigma}=\{\bar s_1,...,\bar s_{n+p+1}\}$ the knot vector obtained by inserting a new knot $\bar s$, $s_i\le \bar s<s_{i+1}$.
Let ${\bf \Omega_u}=\{u_1(s),...,u_{n+p-1}(s)\}$, ${\bf \Omega_v}=\{v_1(s),...,v_{n+p-1}(s)\}$ and ${\bf \bar\Omega_u}=\{\bar u_1(s),...,\bar u_{n+p}(s)\}$, ${\bf \bar\Omega_v}=\{\bar v_1(s),...,\bar v_{n+p}(s)\}$ be the corresponding vectors of functions, where 
\begin{align}
&\bar u_j(s)=u_j(s) \quad\hbox{and}\quad \bar v_j(s)=v_j(s) \quad\hbox{if}\quad j\le i  \notag\\
&\bar u_j(s)=u_{j-1}(s) \quad\hbox{and}\quad \bar v_j(s)=v_{j-1}(s)\quad \hbox{if}\quad j>i. \label{aaa6}%\label{29ot131}
\end{align}
If we denote by $N_{i}^{(p)}(s)$ and $\bar N_{i}^{(p)}(s)$ the GB-splines of order $p$, respectively before and after the knot insertion, and by $r+1$ the multiplicity of $\bar s$ in ${\bf \bar \Sigma}$, then we obtain
\begin{equation}
N_{j}^{(p)}(s)=\alpha_{j,p}\bar N_{j}^{(p)}(s)+\beta_{j+1,p}\bar N_{j+1}^{(p)}(s), \label{a5}
\end{equation}
with, for $p>2$,
\begin{align*}
&\alpha_{j,p}=\begin{cases}
1,&\text{$j\le i-p$,} \\ 
\frac{\displaystyle \delta^{(p-1)}_{j}}{\displaystyle \bar \delta^{(p-1)}_{j}}\alpha_{j,p-1},&\text{$i-p<j<i-r+1$,} \\ 
0,&\text{$j\ge i-r+1$}
\end{cases},
\\
&\beta_{j,p}=\begin{cases}
0,&\text{$j\le i-p+1$,} \\ 
\frac{\displaystyle \delta^{(p-1)}_{j}}{\displaystyle \bar \delta^{(p-1)}_{j+1}}\beta_{j+1,p-1},&\text{$i-p+1<j<i-r+2$,} \\ 
1,&\text{$j\ge i-r+2$}
\end{cases},
\end{align*}
and, for $p=2$,
\begin{align*}
&\alpha_{j,2}=\begin{cases}
1,&\text{$j<i$,} \\ 
\frac{\displaystyle V_{i}(\bar s)}{\displaystyle V_i(s_{i+1})},&\text{$j=i$,} \\ 
0,&\text{$j\ge i+1$}
\end{cases},\\
&\beta_{j,2}=\begin{cases}
0,&\text{$j<i$,} \\ 
\frac{\displaystyle U_{i}(\bar s)}{\displaystyle U_{i}(s_i)},&\text{$j=i$,} \\ 
1,&\text{$j\ge i+1$}
\end{cases},
\end{align*}
where $\delta_{j}^{(p-1)}$ and $\bar\delta_{j}^{(p-1)}$ are the constants defined by \eqref{a3} for ${\bf \Sigma}$ and ${\bf \bar \Sigma}$ respectively, and $U_i(s)$ and $V_i(s)$, $\bar U_i(s)$ and $\bar V_i(s)$, $\bar U_{i+1}(s)$ and $\bar V_{i+1}(s)$ are the generating functions associated to $[s_i,s_{i+1}]$, $[s_i,\bar s]$, $[\bar s,s_{i+1}]$, respectively, and such that $\bar{V}_i(s)=V_i(s), \bar{U}_{i+1}(s)=U_i(s)$.
\end{theorem}

\smallskip
%\noindent
{\bf Proof.} The result can be proved either by generalizing the knot insertion formula found for trigonometric GB-splines  \cite{bracco12b} or by using the properties of piecewise Extended Chebyshev splines (see \cite{mazure04}, \cite{mazure11a} and \cite{mazure11b}). \hfill $\square$

\smallskip
\noindent
Note that knot insertion formulae for a particular case of GB-splines and for slightly different splines can be also found in \cite{wang04} and \cite{wang08}, respectively.

\section{Generalized T-splines}

\medskip
In order to define the GT-splines, we need to briefly recall some definitions and notations about the T-meshes, which are the same used for the classical polynomial case (see, e.g., \cite{cho12} and \cite{daveiga}) and for the TGT-splines in \cite{bracco12b}.

\medskip
%\noindent
Let ${\bf \Sigma^s}=\{s_{-\lfloor p/2\rfloor+1},...,s_{\mu+\lfloor p/2\rfloor}\}$ and ${\bf \Sigma^t}=\{t_{-\lfloor q/2\rfloor+1},...,t_{\nu+\lfloor q/2\rfloor}\}$ be two index vectors, where $\mu,\nu\in \ZZ$, $p,q\in \ZZ$ are equal to or greater than $2$ and, for any real number $k$, $\lfloor k\rfloor$ is the largest integer smaller than or equal to $k$. Analogously, ${\bf \Omega_u^s}$, ${\bf \Omega_v^s}$, ${\bf \Omega_u^t}$ and ${\bf \Omega_v^t}$ are the associated vectors of functions.

\noindent
An {\sl index T-mesh} $M$ is a rectangular partition of the index domain $[-\lfloor p/2\rfloor+1,\mu+\lfloor p/2\rfloor]\times[-\lfloor q/2\rfloor+1,\nu+\lfloor q/2\rfloor]$ such that the vertices have integer coordinates (see Figure 3(a)). In other words, $M$ is the collection of all the elements of such partition, which are called {\sl cells}. Note that, since the elements are rectangular, T-junctions are allowed but L-junctions or I-junctions are not. We call {\sl edge} any segment, either horizontal or vertical, linking two vertices of the mesh. We denote the set of vertices by ${\cal V}$ and by $hE$, $vE$ and $E$ the sets containing only horizontal, only vertical and all the edges respectively. The valence of a vertex $P$ is the number of edges $e\in E$ such that $P\in \partial e$. Finally, we denote by $S$ the union of all the edges and vertices.

\noindent
We define the {\sl active region} $AR_{p,q}$ and {\sl frame region} $FR_{p,q}$ (see Figure 3(b)) as
\begin{equation*}
AR_{p,q}=[1,\mu]\times[1,\nu],
\end{equation*}
and
\begin{align*}
FR_{p,q}=&\bigg(\big[-\lfloor p/2\rfloor+1,1\big]\cup\big[\mu,\mu+\lfloor p/2\rfloor\big]\bigg)\times\big[-\lfloor q/2\rfloor+1,\nu+\lfloor q/2\rfloor\big]\\
&\cup\big[-\lfloor p/2\rfloor+1,\mu+\lfloor p/2\rfloor\big]\times\bigg(\big[-\lfloor q/2\rfloor+1,1\big]\cup\big[\nu,\nu+\lfloor q/2\rfloor \big]\bigg).
\end{align*}

\begin{figure}[H] 
\begin{minipage}[H]{.45\textwidth}
\centering
\includegraphics[scale=0.44]{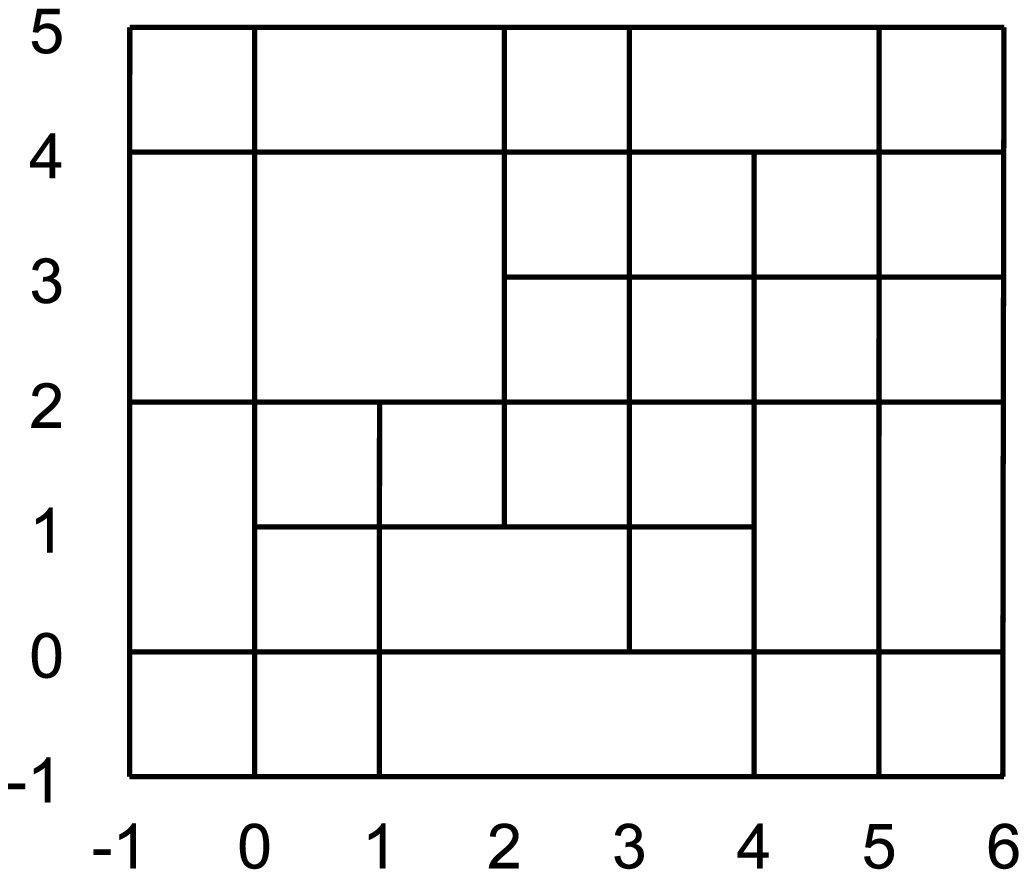}
\caption{\sl (a)}
\end{minipage} 
\hspace{.05\textwidth}
\begin{minipage}[H]{.45\textwidth}
\centering
\includegraphics[scale=0.44]{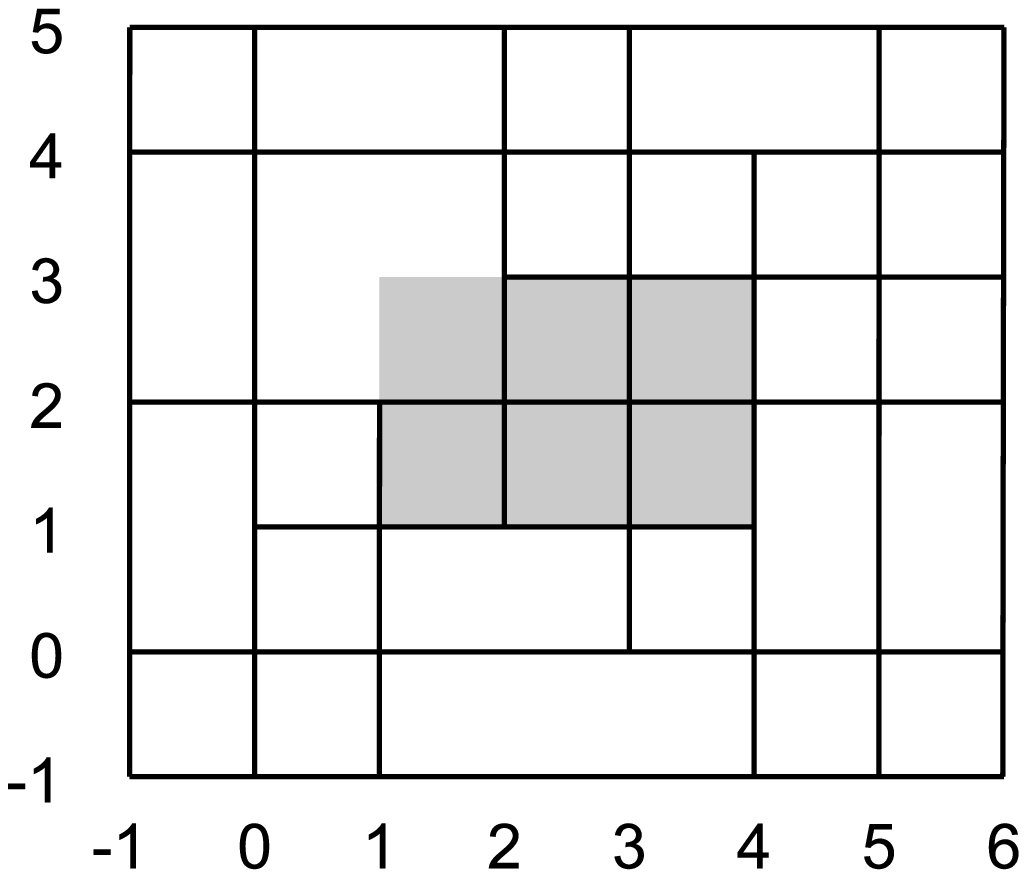}
\caption{\sl (b)}
\end{minipage}
\caption{Figure 3: {\sl A T-mesh (a) in the case $p=q=4$, $\mu=4$ and $\nu=3$ , and (b) the corresponding active region highlighted in gray, with the remaining part representing the frame region.}}
\end{figure}

\begin{definition}\label{B1}
A T-mesh $M$ is admissible for the bi-order $(p,q)$ if $S\cap FR_{p,q}$ includes the segments
\begin{align*}
\{l\}\times\big[-\lfloor q/2\rfloor+1,\nu+\lfloor q/2\rfloor\big]\qquad &\hbox{for}\,\, l=-\lfloor p/2\rfloor+1,...,1,\\
&\hbox{and}\,\, l=\mu,...,\mu+\lfloor p/2\rfloor,\\
\big[-\lfloor p/2\rfloor+1,\mu+\lfloor p/2\rfloor\big]\times\{l\}\qquad &\hbox{for}\,\, l=-\lfloor q/2\rfloor+1,...,1,\\
&\hbox{and}\,\, l=\nu,...,\nu+\lfloor q/2\rfloor,\\
\end{align*}
and all vertices belonging to $(-\lfloor p/2\rfloor+1,\mu+\lfloor p/2\rfloor)\times(-\lfloor q/2\rfloor+1,\nu+\lfloor q/2\rfloor)\cap FR_{p,q}$ have valence 4. $AD_{p,q}$ will denote the set of admissible T-meshes for the bi-order $(p,q)$.
\end{definition}

\begin{definition}\label{B1plus}
A T-mesh $M\in AD_{p,q}$ belongs to $AD_{p,q}^{+}$ if, for any couple of vertices $P_1=(i_1,j_1),P_2=(i_2,j_2)\in {\cal V}$ both belonging to the boundary of a cell and such that $i_1=i_2$ ($j_1=j_2$, resp.), the segment ${i_1}\times (j_1,j_2)$ (${j_1}\times (i_1,i_2)$, resp.) belongs to $S$.
\end{definition}
\noindent
In other words, A T-mesh satisfying the definition \ref{B1plus} does not have any \lq\lq facing\rq\rq T-junctions. While considering this additional requirement is not necessary now, we will need it later to guarantee the equivalence between analysis-suitable and dual-compatible T-meshes (see \cite{daveiga}). 

\noindent
The so-called {\sl anchors}, which are basic to the construction of T-splines, are defined as follows.

\begin{definition}\label{B2}
Given T-mesh $M\in AD_{p,q}$, the set of anchors ${\cal A}_{p,q}(M)$ is defined in the following way:
\begin{itemize}
\item if both $p$ and $q$ are even, ${\cal A}_{p,q}(M)=\{A\in {\cal V}: A\subset AR_{p,q}\}$;
\item if $p$ is odd and $q$ is even, ${\cal A}_{p,q}(M)=\{A\in hE: A\subset AR_{p,q}\}$;
\item if $p$ is even and $q$ is odd, ${\cal A}_{p,q}(M)=\{A\in vE: A\subset AR_{p,q}\}$;
\item if both $p$ and $q$ are odd, ${\cal A}_{p,q}(M)=\{A\in M: A\subset AR_{p,q}\}$.
\end{itemize}
\end{definition}

To each anchor $A\in{\cal A}_{p,q}(M)$ we associate a {\sl global horizontal (vertical) index vector} ${\bf I^s}(A)$ (${\bf I^t}(A)$) and a {\sl local horizontal (vertical) index vector} ${\bf I_l^s}(A)$ (${\bf I_l^t}(A)$), which is a subset of ${\bf I^s}(A)$ (${\bf I^t}(A)$). The construction of these vectors, which depends on the local topology of the T-mesh, is fundamental in the theory of T-splines, and a formal presentation can be found for example in \cite{daveiga}.

\smallskip
%\noindent
The {\sl T-mesh in parameter space} is naturally defined as the partition of the domain $[s_{-\lfloor p/2\rfloor+1},s_{\mu+\lfloor p/2\rfloor}]\times [t_{-\lfloor q/2\rfloor+1},t_{\nu+\lfloor q/2\rfloor}]$ obtained by considering the elements of the form
\begin{equation*}
(s_{i_1},s_{i_2})\times(t_{j_1},t_{j_2})\ne \emptyset, 
\end{equation*}
where $(i_1,i_2)\times(j_1,j_2)\in M$. Let us introduce the notation
\begin{align*}
&{\bf \Sigma^s}({\bf I^s})=\{s_{i}\in {\bf \Sigma^s}:i\in {\bf I^s}\},\\
&{\bf \Sigma^t}({\bf I^t})=\{t_{j}\in {\bf \Sigma^t}:j\in {\bf I^t}\},\\ 
&{\bf \Omega^s_u}({\bf I^s})=\{u^s_{i}\in {\bf \Omega^s_u}:i\in {\bf I^s}\backslash \{i_{h_{p+1}}\}\},\\ 
&{\bf \Omega^t_u}({\bf I^t})=\{u^t_{j}\in {\bf \Omega^t_u}:j\in {\bf I^t}\backslash \{j_{k_{q+1}}\},\\
&{\bf \Omega^s_v}({\bf I^s})=\{v^s_{i}\in {\bf \Omega^s_v}:i\in {\bf I^s}\backslash \{i_{h_{p+1}}\},\\
&{\bf \Omega^t_v}({\bf I^t})=\{v^t_{j}\in {\bf \Omega^t_v}:j\in {\bf I^t}\backslash \{j_{k_{q+1}}\}.
\end{align*}
for any index vectors ${\bf I^s}=\{i_{h_1},...,i_{h_{p+1}}\}\subset\{1,2,...,\mu\}$, ${\bf I^t}=\{j_{k_1},...,j_{k_{q+1}}\}\subset\{1,2,...,\nu\}$. In this way, the global and local index vectors associated to each anchor naturally define corresponding global and local knot and functions vectors.  

\smallskip
\noindent
Then, we define, for each anchor, a bivariate {\sl Generalized T-spline} (GT-spline): 
\begin{align}
N_A(s,t)=&N\big[{\bf \Sigma^s}({\bf I^s_l}(A)),{\bf \Omega_{u}^s}({\bf I^s_l}(A)), {\bf \Omega_{v}^s}({\bf I^s_l}(A))\big](s)\notag \\
&\times N\big[{\bf \Sigma^t}({\bf I^t_l}(A)),{\bf \Omega_{u}^t}({\bf I^t_l}(A)), {\bf \Omega_{v}^t}({\bf I^t_l}(A))\big](t). \label{bb1}
\end{align}
where $N\big[{\bf \Sigma^s}({\bf I^s_l}(A)),{\bf \Omega_{u}^s}({\bf I^s_l}(A)), {\bf \Omega_{v}^s}({\bf I^s_l}(A))\big](s)$ and $N\big[{\bf \Sigma^t}({\bf I^t_l}(A)),{\bf \Omega_{u}^t}({\bf I^t_l}(A)), {\bf \Omega_{v}^t}({\bf I^t_l}(A))\big](t)$ are the univariate GB-splines in the variables $s$ and $t$ constructed on the horizontal and vertical local knot and function vectors associated to $A$. Of course, the TGT-splines introduced in \cite{bracco12b} are a particular case of the just defined GT-splines, obtained by setting $u_i^s(s)=\cos(\omega_i^s(s)), v_i^s(s)=\sin(\omega_i^s(s))$ and $u_j^t(t)=\cos(\omega_j^t(t)), v_j^t(t)=\sin(\omega_j^t(t))$, where $\omega_i^s$ and $\omega_j^t$ are frequencies such that $0<\omega_i^s< \pi/(s_{i+1}-s_i)$ and $0<\omega_j^t< \pi/(t_{j+1}-t_j)$ for any $i,j$.

\smallskip
\noindent
Several properties holding for the polynomial T-splines are also satisfied by the GT-splines.
\begin{property}\label{C0a}
The GT-splines enjoy the following properties, as direct consequence of their definition.
\begin{enumerate}
\item {\bf Continuity:} each blending function $N_A(s,t)$, for any $A\in {\cal A}_{p,q}(M)$, is $(p-m^s_i-1)$ times continuously differentiable with respect to $s$ and $(q-m^t_j-1)$ times continuously differentiable with respect to $t$ at the point $(s_i,t_j)$, where $m^s_i$ and $m^t_j$ are the multiplicities of $s_i$ and $t_j$ in the knot vectors ${\bf \Sigma^s}({\bf I^s}(A))$ and ${\bf \Sigma^t}({\bf I^t}(A))$, respectively.
\item {\bf Positivity:} $N_A(s,t)\ge 0$ for $(s,t)\in\RR^2$, $A\in {\cal A}_{p,q}(M)$ and $p,q\in\NN$, $p,q\ge 2$.
\item {\bf Local support:} if $(s,t)\notin [\min {\bf \Sigma^s}({\bf I^s_l}(A)),\max {\bf \Sigma^s}({\bf I^s_l}(A))]\times[\min {\bf \Sigma^t}({\bf I^t_l}(A)),\max {\bf \Sigma^t}({\bf I^t_l}(A))]$, then $N_A(s,t)= 0$, $A\in {\cal A}_{p,q}(M)$ and $p,q\in\NN$, $p,q\ge 2$.
\item {\bf Linear independence for tensor-product case:} If $M$ is a {\sl tensor-product mesh}, that is, all the vertices have valence 4, then the corresponding blending functions are linearly independent.
\item {\bf Partition of unity for tensor-product case:} If $M$ is a {\sl tensor-product mesh}, then the corresponding blending functions form a partition of unity.
\end{enumerate}
\end{property}

\smallskip
%\noindent
We can use the blending functions \eqref{bb1} to construct a spline surface in the same way as in the polynomial case:
\begin{equation}
\mathbf{T}(s,t)= \frac{\displaystyle \sum_{A\in {\cal A}_{p,q}(M)} \mathbf{T}_A w_A N_A(s,t)}{\displaystyle \sum_{A\in {\cal A}_{p,q}(M)} w_A N_A(s,t)},\label{bb2}
\end{equation}
where $\mathbf{T}_A\in \RR^3$ are given control points and $w_A\in \RR^+$ are the weights. 

\smallskip
Note that the constant function $1$ may not belong to $\hbox{span}\{N_A(s,t) : \,\,A \in \mathcal{A}_{p,q} (M) \}$, and then considering the rational form \eqref{bb2} allows to get the partition of unity property. It can be easily shown that the concepts of standard and semi-standard T-splines (see \cite{sederberg03} and \cite{sederberg04}) can be extended to this non-polynomial setting; therefore, if we construct standard or semi-standard GT-splines, we can avoid using the rational form, which allows us to combine the features of the T-spline approach and the reproduction properties of the GB-splines.

\smallskip
\noindent
As in the GB-splines case, the use of GT-splines is particularly relevant to exactly represent certain shapes, which cannot be obtained with classical T-splines. For example, helical-shaped domains such as helicoids, helicoidal springs and screws can be exactly reproduced by trigonometric GT-splines.

\noindent{\bf Example 1.} Let us consider the helicoid section parametrized by
\begin{align*}
&x(s,t)=s \cos (\omega t)\\
&y(s,t)=s \sin (\omega t)\\
&z(s,t)=t , \qquad r_1\le s\le r_2, \quad 0\le t \le h.
\end{align*}
The helicoid section in Figure 4(a), where $r_1=1/2$, $r_2=1$, $h=6$ and $\omega=3$, is exactly modeled by using suitable GT-splines of bi-order (4,4), which span a space containing $\langle 1,s,\cos(3s), \sin (3s) \rangle \otimes \langle1,t,\cos(3t), \sin (3t)\rangle$.

\begin{figure}[H] 
\begin{minipage}[H]{.45\textwidth}
\centering
\includegraphics[trim=50mm 0mm 50mm 0mm,clip,scale=0.70]{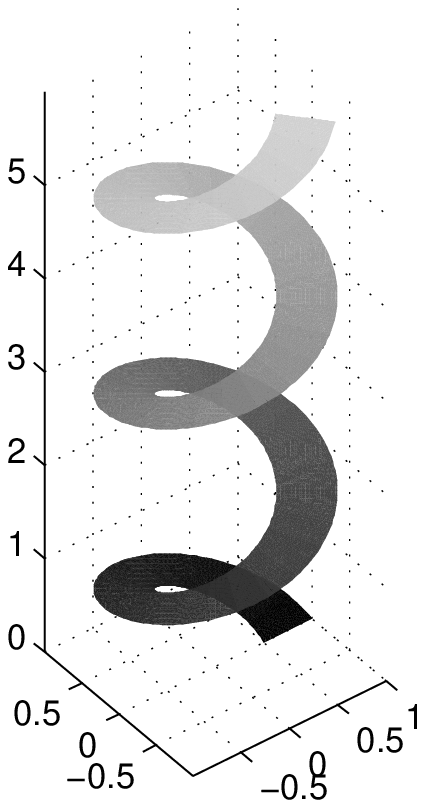}
\caption{\sl (a)}
\end{minipage} 
\hspace{.05\textwidth}
\begin{minipage}[H]{.45\textwidth}
\centering
\includegraphics[trim=20mm 0mm 20mm 0mm,clip,scale=0.70]{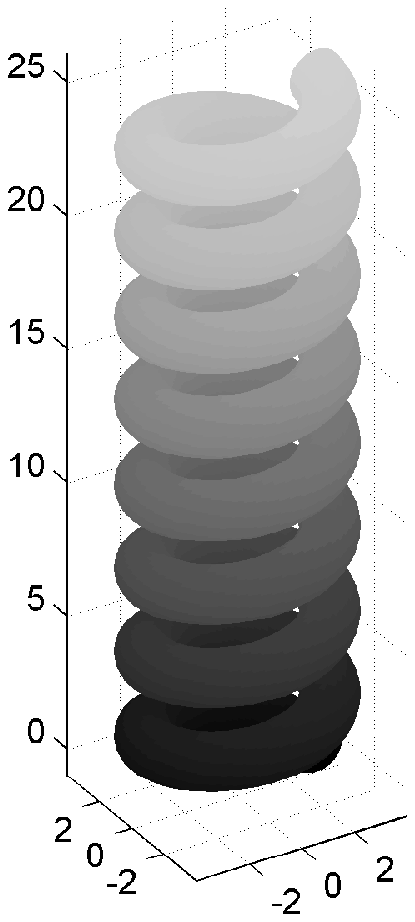}
\caption{\sl (b)}
\end{minipage}
\caption{Figure 4: {\sl A helicoid section (a) in the case $r_1=1/2$, $r_2=1$, $h=10$ and $\omega=3$, and (b) a helicoidal spring in the case $R=3$, $r=1$, $h=8 \pi$, $\omega_s=1$ and $\omega_t=2$.}}
\end{figure}

\noindent{\bf Example 2.} Let us consider a helicoidal spring parametrized by
\begin{align*}
&x(s,t)=(R+r\cos(\omega_s s)) \cos (\omega_t t)\\
&y(s,t)=(R+r\cos(\omega_s s)) \sin (\omega_t t)\\
&z(s,t)=r\sin(\omega_s s)+t, \qquad 0\le s\le 2\pi, \quad 0\le t \le h,
\end{align*}
with $R,r>0$. The helicoidal spring in Figure 4(b), where $R=3$, $r=1$, $h=8 \pi$, $\omega_s=1$ and $\omega_t=2$, is exactly modeled by using suitable GT-splines of bi-order (4,4), which span a space containing $\langle 1,s,\cos(s), \sin (s) \rangle \otimes \langle 1,t,\cos(2t), \sin (2t)\rangle$.

Such shapes are used for example in boundary layer problems on helical-shaped domains (see, e.g., \cite{mostafa} and \cite{rodr}), in modeling of dental implants (see, e.g., \cite{cehreli} and \cite{vstaden}) and in finite element methods on helical-spring models (see, e.g., \cite{jiang} and \cite{toi}). Note that using VMCR T-meshes defined in Section 4.2 guarantees the linear independence of the GT-splines and then makes them suitable to numerically solve the above mentioned problems.

\section{Linear independence of the GT-splines and VMCR T-meshes} 

\subsection{GT-splines and tensor-product splines}

\medskip
The linear independence of the T-splines is a key point for at least one of their main applications, that is, isogeometric analysis (see, e.g., \cite{bazilevs10}). Therefore, the study of linear independence is basic for the theory of the just introduced GT-splines as well, which is the reason why we devote Section 4 to this topic. First we will generalize the results for TGT-splines in \cite{bracco12b}.

\smallskip
\noindent
In general, the situation about linear independence of GT-splines may not coincide with the one of the classical polynomial T-splines. For instance, it has been shown that there are examples where the arguments used to prove the linear dependence of the T-splines do not hold in the case of the TGT-splines (see \cite{bracco12b}). The linear independence of the GT-splines can be studied by examining the relation between them and the tensor-product spline functions associated to the so-called underlying tensor product mesh. Note that the tensor-product GB-spline functions are linearly independent (see Property \ref{C0a}).

\smallskip
\noindent
\begin{definition}\label{C01} 
Given a T-mesh $M$, its {\sl underlying tensor-product mesh} $\hat M$ is the T-mesh with the same index domain of $M$ and obtained by adding to $M$ vertices and edges such that all the vertices have valence $4$ and the knot and functions vectors ${\bf \Sigma^s}$, ${\bf \Sigma^t}$, ${\bf \Omega^s_u}$, ${\bf \Omega^s_v}$, ${\bf \Omega^t_u}$ and ${\bf \Omega^t_v}$ are unvaried.
\end{definition}

%\noindent
Then, we can get the underlying tensor-product mesh $\hat M$ of a T-mesh $M$ by adding edges. Therefore, there is a linear relation between the two sets of GT-splines associated to $M$ and $\hat M$, since adding the edges needed to get $\hat M$ from $M$ corresponds to inserting knots belonging to ${\bf \Sigma^s}$ (${\bf \Sigma^t}$ respectively) in the global knot vectors ${\bf \Sigma^s}({\bf I^s}(A))$ (${\bf \Sigma^t}({\bf I^t}(A))$ respectively) and the corresponding elements in the vectors of functions ${\bf \Omega^s_u}({\bf I^s}(A))$ and ${\bf \Omega^s_v}({\bf I^s}(A))$ (${\bf \Omega^t_u}({\bf I^t}(A))$ and ${\bf \Omega^t_v}({\bf I^t}(A))$ respectively) for $A\in{\cal A}_{p,q}(M)$, which can be handled by using the knot insertion formula. Then these knot insertions must satisfy the requirements of Theorem \ref{00a}. More precisely, condition \eqref{aaa6} must be satisfied when we insert new functions belonging to ${\bf \Omega^{s}_u}$ and ${\bf \Omega^{s}_v}$ (${\bf \Omega^{t}_u}$ and ${\bf \Omega^{t}_v}$ respectively), in the global function vectors ${\bf \Omega^s_u}({\bf I^s}(A))$ and ${\bf \Omega^s_v}({\bf I^s}(A))$ (${\bf \Omega^t_u}({\bf I^t}(A))$ and ${\bf \Omega^t_v}({\bf I^t}(A))$ respectively) of an anchor $A$. By repeatedly applying the knot insertion formula \eqref{a5}, we get a relation of type
\begin{equation}
N_A(s,t)=\sum_{B\in{\cal A}_{p,q}(\hat M)} c_{A,B}\hat N_B(s,t), \qquad A\in {\cal A}_{p,q}(M), \label{ca1}
\end{equation} 
where $\{N_A(s,t)\}_{A\in{\cal A}_{p,q}(M)}$ and $\{\hat N_B(s,t)\}_{B\in{\cal A}_{p,q}(\hat M)}$ are the sets of GT-splines associated to $M$ and $\hat M$, respectively.
If we denote the sets of the anchors of $M$ and $\hat M$ by ${\cal A}_{p,q}(M)=\{A_1,..,A_n\}$ and ${\cal A}_{p,q}(\hat M)=\{\hat A_1,..,\hat A_{\hat n}\}$, \eqref{ca1} can be also written in the matrix form
\begin{equation}
{\bf N}={\bf C}{\bf \hat N}, \label{ca2}
\end{equation}
where ${\bf N}=\big[N_{A_1}(s,t),...,N_{A_n}(s,t)\big]^T$, ${\bf\hat N}=\big[\hat N_{\hat A_1}(s,t),...,\hat N_{\hat A_{\hat n}}(s,t)\big]^T$, and ${\bf C}$ is an $n\times \hat n$ matrix ${\bf C}=(c_{ij})_{i=1,...,n,j=1,...,{\hat n}}$, whose elements are obtained by re-labeling the coefficients $c_{A,B}$ in \eqref{ca1}.
The linear independence of the GT-spline blending functions is equivalent to ${\bf C}$ being a full-rank matrix.

\begin{theorem}\label{C1}
A necessary and sufficient condition for the GT-spline blending functions $\{N_{A_i}(s,t)\}_{i=1}^n$ to be linearly independent is that ${\bf C}$ is full rank.
\end{theorem}

\noindent
{\bf Proof.} See the analogous Theorem in \cite{bracco12b}. \hfill$\square$
 
\smallskip
%\noindent
Given the same T-mesh $M$, the same knot vectors ${\bf \Sigma^s}$ and ${\bf \Sigma^t}$, and as a consequence the same anchors and the same global and local knot vectors, we denote by $\{P_{A_i}(s,t)\}_{i=1,...,n}$ the polynomial T-spline blending functions of bi-degree $(p-1,q-1)$, and by $\{\hat P_{\hat A_i}(s,t)\}_{i=1}^{\hat n}$ the tensor-product B-spline functions associated to the underlying tensor-product mesh $\hat M$. We can obtain also in this case, by repeatedly applying Boehm's knot insertion formula for the polynomial splines (see, e.g., \cite{boehm}), the relation
\begin{equation}
{\bf P}={\bf D}{\bf\hat P}, \label{ca3}
\end{equation}
where ${\bf P}=\big[P_{A_1}(s,t),...,P_{A_n}(s,t)\big]^T$, ${\bf \hat P}=\big[\hat P_{\hat A_1}(s,t),...,\hat P_{\hat A_{\hat n}}(s,t)\big]^T$ and ${\bf D}$ is an $n\times \hat n$ matrix.

\begin{theorem}\label{C2}
A necessary and sufficient condition for the T-spline blending functions $\{P_{A_i}(s,t)\}_{i=1}^n$ to be linearly independent is that ${\bf D}$ is full rank.
\end{theorem}

\noindent
{\bf Proof.} See, e.g., \cite{li12}. \hfill$\square$

\medskip
%\noindent
In the nonpolynomial and the polynomial cases, the linear independence of the respective blending functions is equivalent to the respective matrices ${\bf C}$ and ${\bf D}$ being full rank. We will prove that there is a strong connection between the two matrices. Since the elements of the two matrices ${\bf C}$ and ${\bf D}$ are obtained by a repeated application of the respective knot insertion formulae, we need to understand the relation between their knot insertion formulae, stated in the following Lemma.
\begin{lemma}\label{C3a}
Let be given a knot vector ${\bf \Sigma}=\{s_1,...,s_{n+p}\}$ and another knot vector ${\bf \bar \Sigma}=\{\bar s_1,...,\bar s_{n+p+1}\}$ obtained by inserting a new knot $\bar s$ between $s_i$ and $s_{i+1}$. Moreover, let ${\bf \Omega_u}=\{u_1(s),...,u_{n+p-1}(s)\}$, ${\bf \Omega_v}=\{v_1(s),...,v_{n+p-1}(s)\}$ and ${\bf \bar\Omega_u}=\{\bar u_1(s),...,\bar u_{n+p}(s)\}$, ${\bf \bar\Omega_v}=\{\bar v_1(s),...,\bar v_{n+p}(s)\}$ be their respective vectors of functions, where $\bar u_j(s)=u_j(s)$ and $\bar v_j(s)=v_j(s)$ if $j\le i$ or $\bar u_j(s)=u_{j-1}(s)$ and $\bar v_j(s)=v_{j-1}(s)$ if $j>i$. Consider, for any order $p\ge 2$, the knot insertion formulas for the univariate GB-splines and for the univariate polynomial B-splines
\begin{align*}
&N_j(s)=\alpha_{j,p}\bar N_j(s)+\beta_{j+1,p} \bar N_{j+1}(s), \qquad j=1,...,n,\\
&P_j(s)=\gamma_{j,p}\bar P_j(s)+\eta_{j+1,p} \bar P_{j+1}(s), \qquad j=1,...,n,
\end{align*} 
where by $\bar N_j$ and $\bar P_j$ we denote the GB-splines and the B-splines obtained after the knot insertion, and the coefficients are obtained by using, respectively, \eqref{a5} and the classical knot insertion formula (see, e.g., \cite{boehm}). Then, for $j=1,...,n$, we have
\begin{align}
&\alpha_{j,p},\beta_{j,p},\gamma_{j,p},\eta_{j,p}\ge 0,\notag\\
&\alpha_{j,p}=0 \Longleftrightarrow \gamma_{j,p}=0,\label{C02apr1}
\beta_{j,p}=0 \Longleftrightarrow \eta_{j,p}=0.
\end{align}
\end{lemma}

\smallskip
\noindent
{\bf Proof.} The result, analogously to the trigonometric case considered in \cite{bracco12b}, follows from the analysis of the expressions of $\alpha_{j,p},\beta_{j,p},\gamma_{j,p},\eta_{j,p}$, $j=1,...,n$. \hfill $\square$

\medskip
\noindent
This Lemma allows us to establish a connection between the matrices ${\bf C}$ and ${\bf D}$. Combining Lemma \ref{C3a} with the same arguments in the proof of Corollary 4.6 in \cite{bracco12b}, we prove below that the sparsity pattern of the two matrices coincide.

\medskip
\begin{theorem}\label{C3}
Let $M\in AD_{p,q}$, and let $\{N_{A_j}\}_{j=1,...,n}$ and $\{P_{A_j}\}_{j=1,...,n}$ be the sets of the GT-spline blending functions and of the T-spline blending functions associated to $M$, respectively. Moreover, let $\{\hat N_{\hat A_j}\}_{j=1,...,\hat n}$ and $\{\hat P_{\hat A_j}\}_{j=1,...,\hat n}$ be the sets of the GT-spline blending functions and of the T-spline blending functions associated to the underlying tensor-product mesh $\hat M$. If we denote by ${\bf C}$ and ${\bf D}$ the matrices expressing the relation between the functions $\{N_{A_j}\}_{j=1,...,n}$ and $\{\hat N_{\hat A_j}\}_{j=1,...,\hat n}$, and between $\{P_{A_j}\}_{j=1,...,n}$ and $\{\hat P_{\hat A_j}\}_{j=1,...,\hat n}$, defined in \eqref{ca2} and \eqref{ca3}, we have that
\begin{equation}
c_{ij}=0 \Longleftrightarrow d_{ij}=0, \qquad i=1,...,n, j=1,...,\hat n. \label{ca4}
\end{equation}
\end{theorem}

\smallskip
\noindent
{\bf Proof.} The proof is analogous to the one for the trigonometric case considered in \cite{bracco12b}. \hfill $\square$

\subsection{VMCR T-meshes}

\medskip
Using Theorem \ref{C3} and the concept of {\sl column reduction} (employed for bicubic splines in \cite{li12}) we will now define the class of T-meshes which guarantees the linear independence both for the classical polynomial T-splines and for the GT-splines.

\noindent
Let us recall the procedure of column reduction. Given a matrix ${\bf C}$, if all the elements of the $i$-th row are zeros except the $j$-th one, then we call the $j$-th column {\sl innocuous}. The column reduction procedure consists of removing from ${\bf C}$ all the innocuous columns and all the zero rows left after the column removal. The following lemma provides a sufficient condition on the result of column reduction and the rank of the considered matrix.

\begin{lemma}\label{C5}
Given an $m\times n$ matrix ${\bf Q}$ ($m\le n$), if the column reduction procedure applied to ${\bf Q}^T$ gives as result the void matrix, then ${\bf Q}$ is a full rank matrix.
\end{lemma}

\noindent
{\bf Proof.} See, e.g., \cite{li12}.\hfill$\square$

\smallskip
\noindent
If we apply the column reduction to the matrices ${\bf C}^T$ and ${\bf D}^T$ defined in \eqref{ca2} and \eqref{ca3}, we can state the following result.

\begin{proposition}\label{C4}
Let ${\bf C}$ and ${\bf D}$ be the matrices defined in \eqref{ca2} and \eqref{ca3}, and let ${\bf C_{CR}}$ and ${\bf D_{CR}}$ be the matrices obtained by applying the column reduction procedures to ${\bf C}^T$ and ${\bf D}^T$, respectively. Then we have
\begin{equation}
{\bf C_{CR}}=\emptyset \Longleftrightarrow {\bf D_{CR}}=\emptyset, \label{ca5}
\end{equation}
where $\emptyset$ stands for the void $0\times 0$ matrix.
\end{proposition}

%\noindent
{\bf Proof.} The equivalence \eqref{ca5} is a direct consequence of Theorem \ref{C3} and of the definition of column reduction.\hfill$\square$

\smallskip
\noindent
As a consequence, the following Corollary holds.

\begin{corollary}\label{C6}
There exists a class of T-meshes for which both the associated GT-spline blending functions of bi-order $(p,q)$ and the T-spline blending functions of bi-degree $(p-1,q-1)$ are linearly independent. This class is defined as the class of T-meshes such that the matrix ${\bf C_{CR}}$ obtained by applying the column reduction procedure to ${\bf C}^T$ (equivalently, the matrix ${\bf D_{CR}}$ obtained by applying the column reduction procedure to ${\bf D}^T$) is the void matrix. We will call it the class of {\sl VMCR T-meshes (Void Matrix after Column Reduction T-meshes)}. Observe that all the tensor-product meshes belong to this class.
\end{corollary}

%\noindent
{\bf Proof.} Let ${\bf C_{CR}}$ and ${\bf D_{CR}}$ be the matrices obtained by applying the column reduction procedure to ${\bf C}^T$ and ${\bf D}^T$, defined in \eqref{ca2} and \eqref{ca3}, respectively. Let us consider the class of T-meshes such that the matrix ${\bf C_{CR}}$ is the void matrix or, equivalently, such that the matrix ${\bf D_{CR}}$ is the void matrix. In fact, if one of these conditions is satisfied for a T-mesh $M$, by Proposition \ref{C4} also the other is satisfied, and by Lemma \ref{C5} ${\bf C}$ and ${\bf D}$ are full rank. This implies, by Theorems \ref{C1} and \ref{C2}, that the GT-spline and the T-spline blending functions associated to $M$ are linearly independent. \hfill$\square$

\medskip
{\bf Remark.} From the definition of the class, it's clear that in order to check whether or not a T-mesh $M$ is VMCR we do not need to completely compute the corresponding matrix ${\bf C}$ (equivalently ${\bf D}$ for the polynomial case). In fact, the column reduction procedure depends only on the sparsity pattern of the matrix, that is, on which of its elements are zero; as a consequence, we do not need to compute the values of the elements of the matrix, but just to check whether or not they are null. Such information can be obtained only by using the topology of $M$ and the multiplicities of the knots in the local vectors, as we will show now. 

\noindent
Let us define the vectors ${\bf \bar I_l^s}(A)$, ${\bf \bar I_l^t}(A)$ obtained from ${\bf I_l^s}(A)$ and ${\bf I_l^t}(A)$ with the following procedure:
\begin{itemize}
\item first construct ${\bf \bar I_l^s}(A)$ by adding to ${\bf I_l^s}(A)$ the elements $i\in\{1,2,...,\mu\}$ such that $i\notin{\bf I_l^s}(A)$ and $\min({\bf I_l^s}(A))< i < \max({\bf I_l^s}(A))$;
\item remove enough elements at the beginning and at the end of ${\bf \bar I_l^s}(A)$ (starting from the first and from the last one, respectively) so that the multiplicities of the first and last element in ${\bf \Sigma}({\bf \bar I_l^s}(A))$ are the same as in ${\bf \Sigma}({\bf I_l^s}(A))$. 
\end{itemize}
The vector ${\bf \bar I_l^t}(A)$ can be constructed analogously adding elements of $\{1,2,...,\nu\}$ to ${\bf I_l^t}(A)$. We present a basic example in order to clarify this procedure.  

\smallskip
\noindent
{\bf Example 3.} Let us assume $p=4$, $\mu=10$ and 
\begin{equation*}
{\bf \Sigma^s}=\{0,0,0,0,1/8,1/8,1/8,1/4,1/2,3/4,1,1,1,1\},
\end{equation*} 
and let us consider an anchor $A$ such that ${\bf I_l^s}(A)=\{3, 5, 7, 8,9\}$, and then ${\bf \Sigma^s}({\bf I_l^s}(A))=\{1/8, 1/8, 1/2, 3/4,1\}$. In this case, the two steps to get ${\bf \bar I_l^s}(A)$ correspond to:
\begin{itemize}
\item adding to ${\bf I_l^s}(A)$ the element $4$:
\begin{equation*}
{\bf \bar I_l^s}(A)=\{3,4, 5, 6, 7, 8,9\}\Longrightarrow{\bf \Sigma}({\bf \bar I_l^s}(A))=\{1/8, 1/8, 1/8, 1/4, 1/2, 3/4,1\};
\end{equation*}
\item remove one element at the beginning of ${\bf \bar I_l^s}(A)$, since the multiplicity of $1/8$ is different in the vectors
\begin{align*}
&{\bf \Sigma}({\bf \bar I_l^s}(A))=\{1/8, 1/8, 1/8, 1/4, 1/2, 3/4,1\}, \\
&{\bf \Sigma}({\bf I_l^s}(A))=\{1/8, 1/8, 1/2, 3/4,1\}
\end{align*}
$\Longrightarrow$ remove the element $3$ from ${\bf \bar I_l^s}(A)$ and then
\begin{align*}
&{\bf \bar I_l^s}(A)=\{4, 5, 6, 7, 8,9\},\\
&{\bf \Sigma}({\bf \bar I_l^s}(A))=\{1/8, 1/8, 1/4, 1/2, 3/4,1\}.
\end{align*} 
\end{itemize}

\begin{lemma}\label{C17set1}
Let $A\in {\cal A}_{p,q}(M)$; by definition \eqref{bb1}, the associated GT-spline blending function can be represented in the form
\begin{align*}
N_A(s,t)=&N\big[{\bf \Sigma^s}({\bf I^s_l}(A)),{\bf \Omega_{u}^s}({\bf I^s_l}(A)), {\bf \Omega_{v}^s}({\bf I^s_l}(A))\big](s)\notag \\
&\times N\big[{\bf \Sigma^t}({\bf I^t_l}(A)),{\bf \Omega_{u}^t}({\bf I^t_l}(A)), {\bf \Omega_{v}^t}({\bf I^t_l}(A))\big](t).%\label{21giu2}
%&\hbox{or}\notag\\
\end{align*}
We have that
\begin{align*}
&N\big[{\bf \Sigma^s}({\bf I^s_l}(A)),{\bf \Omega_{u}^s}({\bf I^s_l}(A)), {\bf \Omega_{v}^s}({\bf I^s_l}(A))\big](s)=\sum_{{\bf I^s}\in H(A)} a^s_{{\bf I^s}} N\big[{\bf \Sigma^s}({\bf I^s}),{\bf \Omega_{u}^s}({\bf I^s}), {\bf \Omega_{v}^s}({\bf I^s})\big](s), \\
&N\big[{\bf \Sigma^t}({\bf I^t_l}(A)),{\bf \Omega_{u}^t}({\bf I^t_l}(A)), {\bf \Omega_{v}^t}({\bf I^t_l}(A))\big](t)=\sum_{{\bf I^t}\in K(A)} a^t_{{\bf I^t}} N\big[{\bf \Sigma^t}({\bf I^t}),{\bf \Omega_{u}^t}({\bf I^t}), {\bf \Omega_{v}^t}({\bf I^t})\big](t),
\end{align*}
and, as a consequence,
\begin{equation*}
N_A(s,t)=\sum_{{\bf I^s}\in H(A)}\sum_{{\bf I^t}\in K(A)}a^s_{{\bf I^s}}a^t_{{\bf I^t}} N\big[{\bf \Sigma^s}({\bf I^s}),{\bf \Omega_{u}^s}({\bf I^s}), {\bf \Omega_{v}^s}({\bf I^s})\big](s)N\big[{\bf \Sigma^t}({\bf I^t}),{\bf \Omega_{u}^t}({\bf I^t}), {\bf \Omega_{v}^t}({\bf I^t})\big](t),%\label{21giu3}
\end{equation*}
where $a^s_{{\bf I^s}}>0$, $a^t_{{\bf I^t}}>0$, and $H(A)$, $K(A)$ are sets of index vectors defined by
\begin{equation*}
H(A)=\{{\bf I^s}:{\bf I^s}\subset {\bf \bar I^s_l}(A) \}, \quad K(A)=\{{\bf I^t}:{\bf I^t}\subset {\bf \bar I^t_l}(A) \}.
\end{equation*} 
\end{lemma}

\smallskip
\noindent
{\bf Proof.} The lemma is a direct consequence of the knot insertion formula applied to the GT-splines associated to the T-mesh $M$. \hfill$\square$

\smallskip
\noindent
The following Proposition follows immediately from Lemma \ref{C17set1}.

\begin{proposition}\label{1mar14}
For any $A\in {\cal A}_{p,q}(M)$
\begin{equation*}
N_A(s,t)=\sum_{\hat A\in {\cal A}_{p,q}(\hat M)} c_{\hat A}N_{\hat A}(s,t),%\label{17set2}
\end{equation*}
where $c_{\hat A}>0$ if and only if ${\bf I_l^s}(\hat A)\in H(A)$ and ${\bf I_l^t}(\hat A)\in K(A)$, and $H(A)$ and $K(A)$ are defined in Lemma \ref{C17set1}.
\end{proposition}

\medskip
\noindent
Note that ${\bf \bar I_l^s}(A)$ and ${\bf \bar I_l^t}(A)$, and then $H(A)$ and $K(A)$, can be obtained starting only from the knowledge of the local index vectors and of the knot multiplicities in the local vectors, without the knot insertion formula and its related computations. By Proposition \ref{1mar14}, the same information also makes it possible to compute a matrix ${\bf S}$ having the same size and sparsity pattern as ${\bf C}$: for any $A \in {\cal A}_{p,q}(M)$ and $\hat A\in{\cal A}_{p,q}(\hat M)$, the element of ${\bf S}$ in the row corresponding to $A$ and in the column corresponding to $\hat A$ is
\begin{itemize}
\item $1$ if ${\bf I_l^s}(\hat A)\in H(A)$ and ${\bf I_l^t}(\hat A)\in K(A)$;
\item $0$ otherwise. 
\end{itemize}
Once this matrix has been computed, we can check if the T-mesh is VMCR by applying the column reduction to it, since it has the same sparsity pattern of $C$.

%%%%%%%%%%%%%%%%%%%%

\subsection{Weakly dual-compatible T-meshes}
\medskip
%\noindent
We will now study the class of VMCR T-meshes. In particular, we will provide a simple characterization for a noteworthy sub-class of T-meshes, which we will call {\sl weakly dual-compatible T-meshes}. We will prove that the class of {\sl dual-compatible/analysis-suitable} T-meshes, the most known class of T-meshes guaranteeing the linear independence of the associated T-spline blending functions, is included in the class of {weakly dual-compatible T-meshes}. Moreover, we will show that there exist non-analysis-suitable T-meshes in the class of weakly dual-compatible T-meshes.

\smallskip
%\noindent
Let us recall what analysis-suitable means (see also \cite{li12} and \cite{daveiga}). Given a T-mesh $M$ (in the index space), let us consider a T-junction $T=(\bar \i,\bar \j)$ belonging to the active region $AR_{p,q}$ and with valence $3$, and assume it is of type \lq\lq$\dashv$\rq\rq, that is, two opposite vertical edges and one horizontal edge from left intersects at $T$. Moreover, let us consider the set of indices $hJ(\bar \j)$ and let $i_1,...i_{p}$ be the $p$ consecutive indices extracted from $hJ(\bar \j)$ such that $i_k=\bar \i$, with $k=\lceil p/2 \rceil$. Then the {\sl horizontal extension} $hext_{p,q}(T)$, with respect to the bi-order $(p,q)$, is defined as the union of the {\sl face extension} $hext_{p,q}^f(T)$ and of the {\sl edge extension} $hext_{p,q}^e(T)$, which are determined as follows:
\begin{equation*}
hext_{p,q}^f(T)=[\bar \i,i_{p}]\times \bar \j,\qquad hext_{p,q}^e(T)=[i_1,\bar \i]\times \bar \j.
  \end{equation*}
We can define analogously the extensions for the other types of T-junctions (see Figure 5 for an example).

\begin{figure}[h] 
\begin{minipage}[H]{.45\textwidth}
\centering
\includegraphics[scale=0.44]{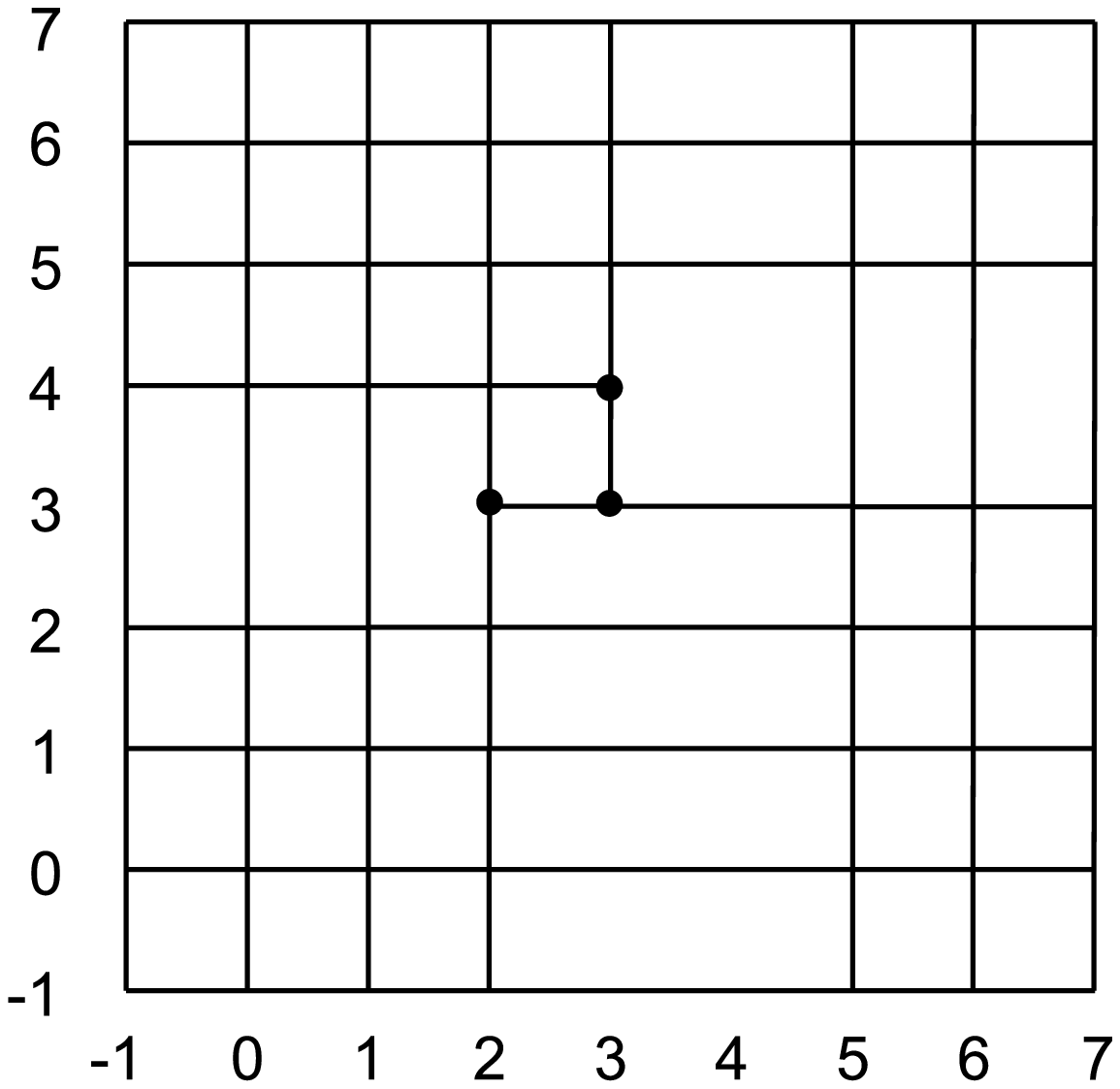}
\caption{\sl (a)}
\end{minipage} 
\hspace{.05\textwidth}
\begin{minipage}[H]{.45\textwidth}
\centering
\includegraphics[scale=0.44]{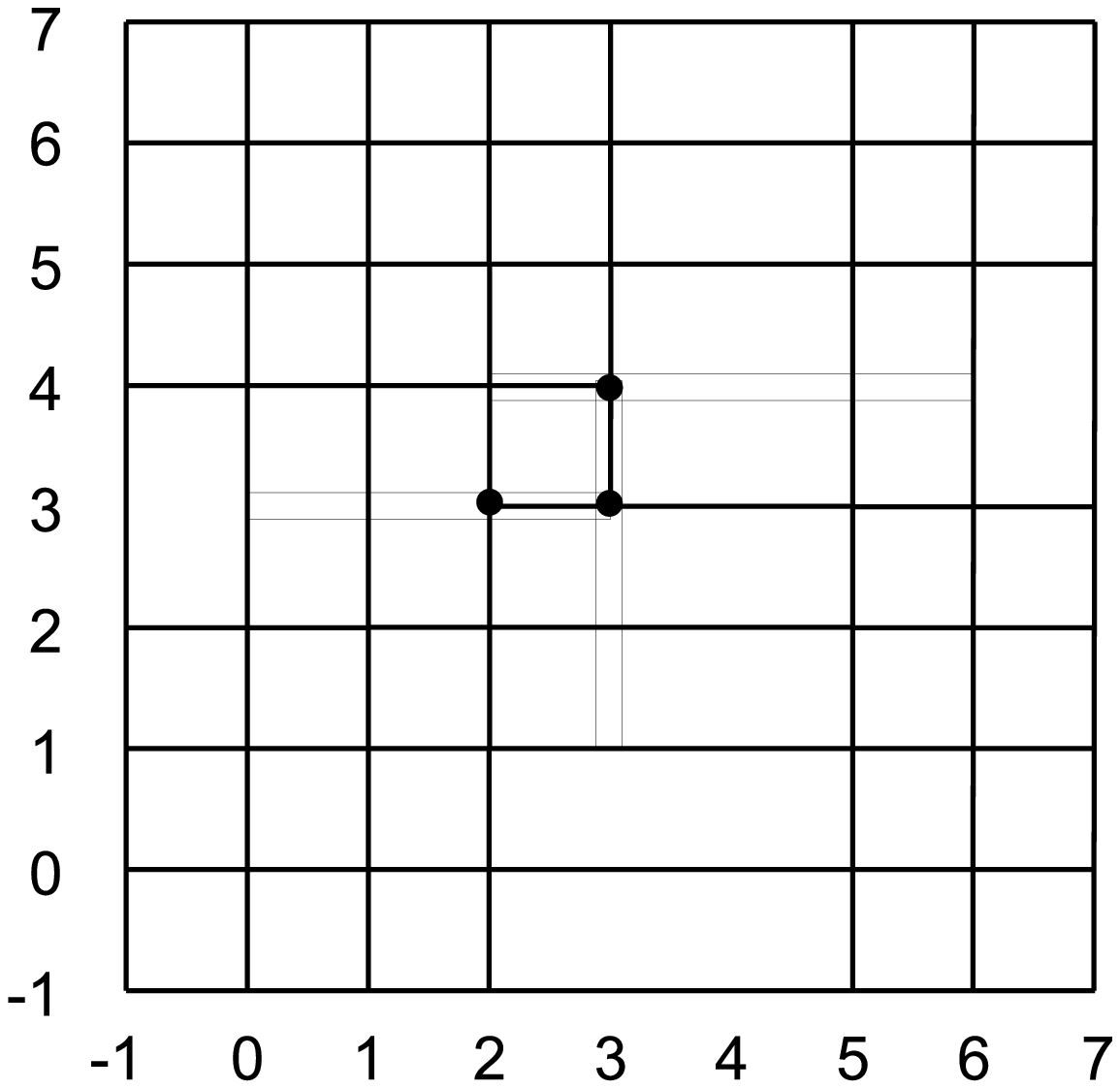}
\caption{\sl (b)}
\end{minipage}
\caption{Figure 5: {\sl A T-mesh (a), with $p=q=4$, and (b) the T-node extensions of its T-nodes.}}
\end{figure}

\begin{definition}\label{C7}
If no horizontal extension with respect to the bi-order $(p,q)$ intersects a vertical extension with respect to the bi-order $(p,q)$, the T-mesh is called {\sl analysis suitable} with respect to the bi-order $(p,q)$ (see Figure 6 for an example).
\end{definition}

\begin{figure}[h] 
\begin{minipage}[H]{.45\textwidth}
\centering
\includegraphics[scale=0.44]{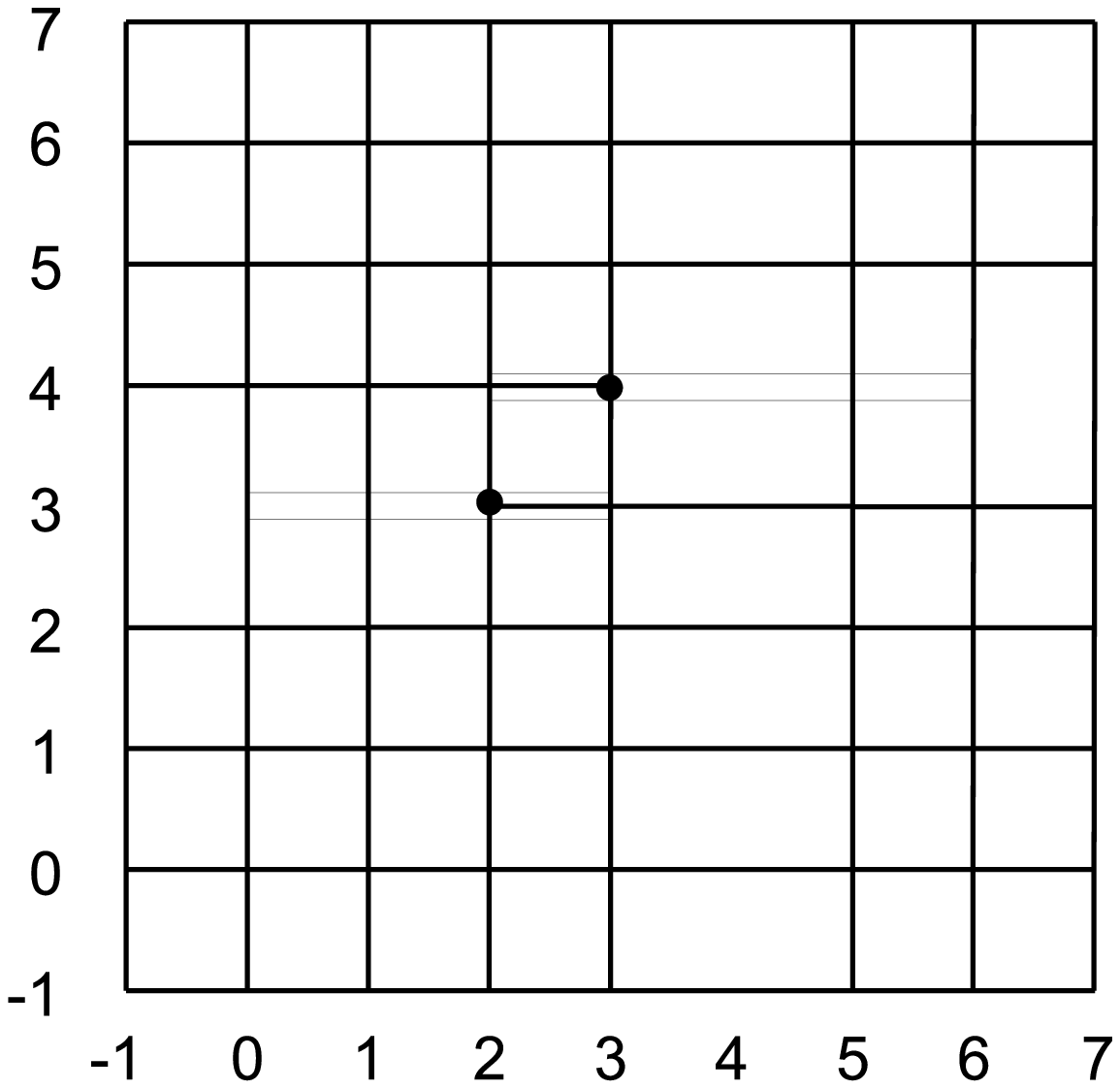}
\caption{\sl (a)}
\end{minipage} 
\hspace{.05\textwidth}
\begin{minipage}[H]{.45\textwidth}
\centering
\includegraphics[scale=0.44]{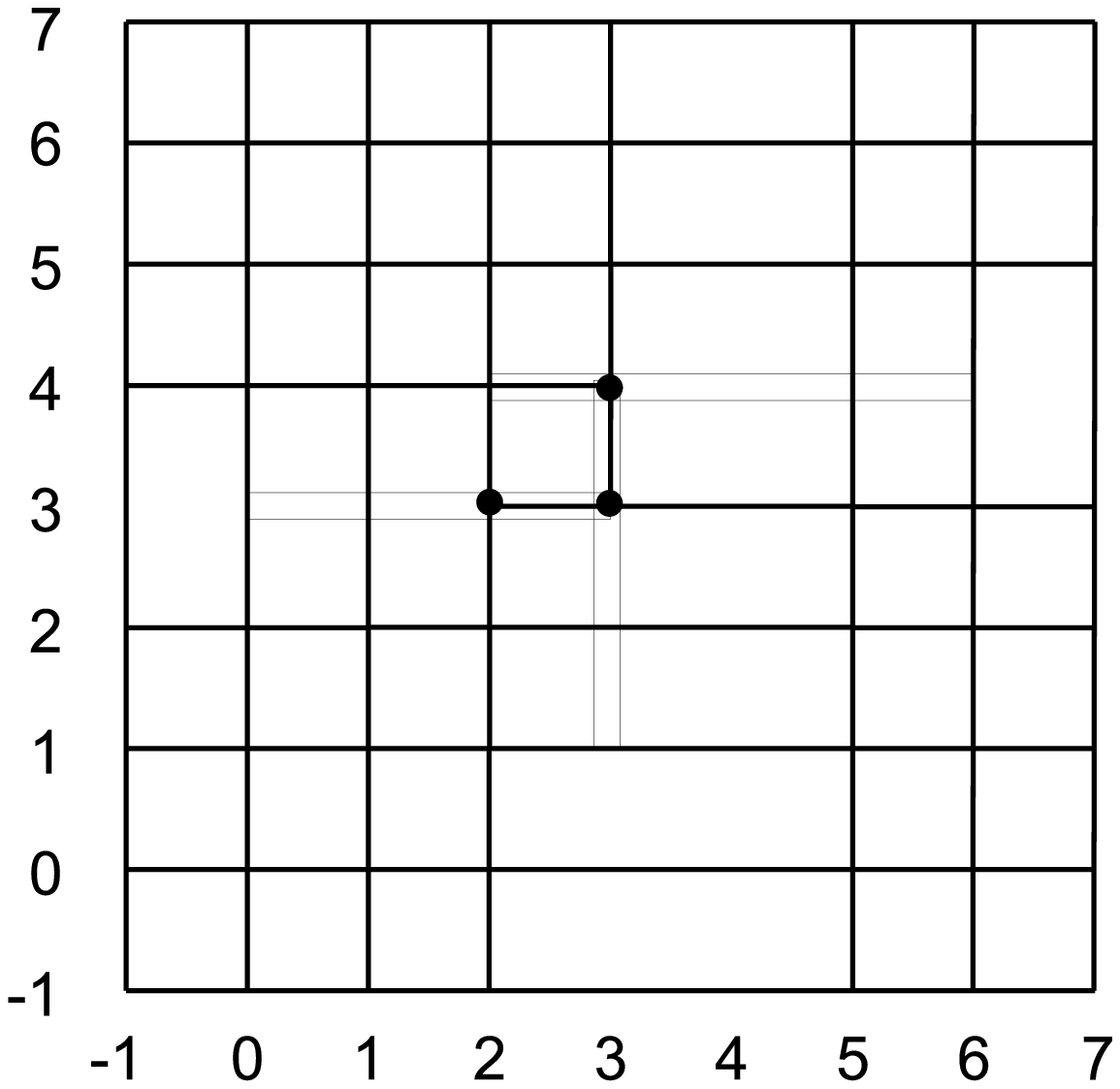}
\caption{\sl (b)}
\end{minipage}
\caption{Figure 6: {\sl An example of (a) analysis-suitable and (b) non-analysis-suitable T-mesh (with $p=q=4$).}}
\end{figure}

\smallskip
The class of analysis-suitable T-meshes coincides with another one: dual-compatible T-meshes. This class was introduced in \cite{cho12} and its equivalence to analysis-suitable T-meshes was proved, for a general bi-degree, in \cite{daveiga}. Since this equivalence will be the key to study the relationship between the analysis-suitable and the weakly dual-compatible T-meshes, let us recall the definition of dual-compatible T-mesh.

\noindent
Let $M\in AD^+_{p,q}$ and let $A_1$ and $A_2$ be two anchors with local horizontal index vectors ${\bf I^s_l}(A_1)=\{i_1^s(A_1),...,i_{p+1}^s(A_1)\}$ and ${\bf I^s_l}(A_2)=\{i_1^s(A_2),...,i_{p+1}^s(A_2)\}$. We say that $A_1$ and $A_2$ {\sl overlap horizontally} if
\begin{align*}
& \forall k\in {\bf I^s_l}(A_1), \,\, i_1^s(A_2)\le k \le i_{p+1}^s(A_2)\Rightarrow k\in {\bf I^s_l}(A_2),\\
&\forall k\in {\bf I^s_l}(A_2), \,\, i_1^s(A_1)\le k \le i_{p+1}(A_1)\Rightarrow k\in {\bf I^s_l}(A_1).%\label{C27mar1}
\end{align*}
Analogously, if ${\bf I^t_l}(A_1)=\{i_1^t(A_1),...,i_{q+1}^t(A_1)\}$ and ${\bf I^t_l}(A_2)=\{i_1^t(A_2),...,i_{q+1}^t(A_2)\}$ are the vertical index vectors of $A_1$ and $A_2$, we we say that $A_1$ and $A_2$ {\sl overlap vertically} if
\begin{align*}
& \forall h\in {\bf I^t_l}(A_1), \,\, i_1^t(A_2)\le h \le i_{p+1}^t(A_2)\Rightarrow h\in {\bf I^t_l}(A_2),\\
&\forall h\in {\bf I^t_l}(A_2), \,\, i_1^t(A_1)\le h \le i_{p+1}^t(A_1)\Rightarrow h\in {\bf I^t_l}(A_1).%\label{C27mar2}
\end{align*}
Moreover, the anchors $A_1$ and $A_2$ are said to {\sl partially overlap} if they overlap either horizontally or vertically.

\smallskip
\begin{definition}\label{C27mar3}
A T-mesh $M\in AD^+_{p,q}$ is dual-compatible with respect to the bi-order $(p,q)$ if any two anchors $A_1,A_2\in {\cal A}_{p,q}(M)$ partially overlap.
\end{definition}

\smallskip
\noindent
Let us now introduce the class of {\sl weakly dual-compatible T-meshes}.

\begin{definition}\label{NN12mag1}
Two anchors $A_1,A_2\in{\cal A}_{p,q}(M)$ are {\sl left-horizontally shifted} ({\sl right-horizontally shifted}, respectively) if conditions \eqref{nnmag121a} (\eqref{nnmag121b}, respectively) are verified:
\begin{align}\label{nnmag121a}
&\vert\min {\bf I_l^s}(A_1)-\min {\bf I_l^s}(A_2)\vert >0\quad \hbox{and, if} \quad \min {\bf \Sigma^s}({\bf I_l^s}(A_1))=\min {\bf \Sigma^s}({\bf I_l^s}(A_2)),\\ \notag
&\hbox{their respective multiplicities in the vectors}\quad{\bf \Sigma^s}({\bf I_l^s}(A_1))\,\,\hbox{and}\,\,\, {\bf \Sigma^s}({\bf I_l^s}(A_2)) \quad \hbox{are different,}
\end{align}
\begin{align}\label{nnmag121b}
&\vert \max {\bf I_l^s}(A_1)-\max {\bf I_l^s}(A_2)\vert >0 \quad \hbox{and, if} \quad \max {\bf \Sigma^s}({\bf I_l^s}(A_1))=\max {\bf \Sigma^s}({\bf I_l^s}(A_2)),\\ \notag
&\hbox{their respective multiplicities in the vectors}\quad{\bf \Sigma^s}({\bf I_l^s}(A_1))\,\,\hbox{and}\,\,\, {\bf \Sigma^s}({\bf I_l^s}(A_2)) \quad \hbox{are different.}
\end{align}
Analogously, two anchors $A_1,A_2\in{\cal A}_{p,q}(M)$ are {\sl down-vertically shifted} ({\sl up-horizontally shifted}, respectively) if conditions \eqref{nnmag121c} (\eqref{nnmag121d}, respectively) are verified:
\begin{align}\label{nnmag121c}
&\vert \min {\bf I_l^t}(A_1)-\min {\bf I_l^t}(A_2)\vert>0 \quad \hbox{and, if} \quad \min {\bf \Sigma^t}({\bf I_l^t}(A_1))=\min {\bf \Sigma^t}({\bf I_l^t}(A_2)),\\ \notag
&\hbox{their respective multiplicities in the vectors}\quad{\bf \Sigma^t}({\bf I_l^t}(A_1))\,\,\hbox{and}\,\,\, {\bf \Sigma^t}({\bf I_l^t}(A_2)) \quad \hbox{are different,}
\end{align}
\begin{align}\label{nnmag121d}
&\vert \max {\bf I_l^t}(A_1)-\max {\bf I_l^t}(A_2)\vert>0 \quad \hbox{and, if} \quad \max {\bf \Sigma^t}({\bf I_l^t}(A_1))=\max {\bf \Sigma^t}({\bf I_l^t}(A_2)),\\ \notag
&\hbox{their respective multiplicities in the vectors}\quad{\bf \Sigma^t}({\bf I_l^t}(A_1))\,\,\hbox{and}\,\,\, {\bf \Sigma^t}({\bf I_l^t}(A_2)) \quad \hbox{are different.}
\end{align}
\end{definition}

\begin{definition}\label{NN12mag2}
A T-mesh is {\sl weakly dual-compatible}:
\begin{itemize}
\item of type RD if any two distinct anchors are either right-horizontally or down-vertically shifted;
\item of type RU if any two distinct anchors are either right-horizontally or up-vertically shifted;
\item of type LD if any two distinct anchors are either left-horizontally or down-vertically shifted;
\item of type LU if any two distinct anchors are either left-horizontally or up-vertically shifted.
\end{itemize}
\end{definition}

\smallskip
\noindent
In order to show that the class of VMCR T-meshes contains the weakly dual-compatible ones for any bi-order $(p,q)$, first we introduce the notion of {\sl influence sub-matrix of a set of anchors ${\cal A}$}, a generalization of the influence sub-graphs used by Li and his co-authors in \cite{li12}. We will give the definitions and the following results referring to the GT-spline blending functions, but it can be easily verified, by using Theorem \ref{C3}, that they hold for the polynomial T-splines as well.

\smallskip
\begin{definition}\label{C28mar1}
Given a set of anchors ${\cal A}\subset {\cal A}_{p,q}(M)$, the influence submatrix of ${\cal A}$ (for the GT-spline blending functions), denoted by ${\bf C}({\cal A})$, is obtained from the matrix ${\bf C}$ defined in \eqref{ca2} by removing the rows corresponding to the anchors not belonging to ${\cal A}$ and all the zero columns left after the rows removal.
\end{definition}

\noindent
We observe that, since each ${\bf C}({\cal A})$ is essentially a submatrix of ${\bf C}$ defined in \eqref{ca2}, each of its rows corresponds to an anchor in ${\cal A}_{p,q}(M)$ and each of its columns corresponds to an anchor in ${\cal A}_{p,q}(\hat M)$, where $\hat{M}$ is the underlying tensor product mesh of $M$.

\begin{definition}\label{C28mar2}
The influence submatrix ${\bf C}({\cal A})$ of a given set of anchors ${\cal A}$ is called a 2-influence submatrix if each column has at least $2$ non-zero elements.
\end{definition}

\begin{lemma}\label{C28mar3}
Given a T-mesh $M\in AD_{p,q}$, if for any set of anchors ${\cal A}\subset {\cal A}_{p,q}(M)$ the influence submatrix ${\bf C}({\cal A})$ for the GT-spline functions is not a 2-influence submatrix, then $M$ belongs to the class of VMCR T-meshes.
\end{lemma}

\smallskip
\noindent
{\bf Proof.} The matrices we obtain at each step of the procedure of column reduction applied to ${\bf C}^T$ can be considered as transpose matrices of influence submatrices, since removing columns and the zero rows left after the columns removal from ${\bf C}^T$ is equivalent to removing rows and the zero columns left after the rows removal from ${\bf C}$. Therefore, since by hypothesis the transpose of each of these matrices is not a 2-influence submatrix, we have that, at each step of the procedure, the obtained matrix has at least a row with no more than one non-zero element. As a consequence, a further column reduction can be always performed, until we reach the void matrix, which proves the Lemma. \hfill $\square$

\medskip
\begin{theorem}\label{NN12mag3}
A weakly dual-compatible T-mesh $M$ is VMCR.
\end{theorem}

\smallskip
\noindent
{\bf Proof.} We prove it for weakly dual-compatible T-meshes of type RD (analogous arguments can be used for the other types). Let us consider any set of anchors ${\cal A}\subset {\cal A}_{p,q}(M)$, and  the set of anchors of $\hat M$ corresponding to the columns of ${\bf C}({\cal A})$, denoted by $\hat {\cal A}$. By Lemma \ref{C28mar3}, in order to prove the Theorem it's sufficient to show that there is at least an anchor $\bar A\in\hat {\cal A}$ such that the corresponding column of ${\bf C}({\cal A})$ has only a non-zero entry. We claim that this anchor can be chosen as the one satisfying, for any $\hat A\in\hat {\cal A}$, one of the two following conditions:
\begin{align}
&i^t_{q+1}(\hat A)>i^t_{q+1}(\bar A),\label{129a}\\
&i^t_{q+1}(\hat A)=i^t_{q+1}(\bar A)\quad \hbox{and} \quad i^s_{p+1}(\hat A)< i^s_{p+1}(\bar A)\label{129c}.
\end{align}
Roughly speaking, these two conditions means that we first choose the lowest anchor in $\hat{\cal A}$ and then, if it is not unique, the rightmost one. The ability to choose $\bar A$ is guaranteed by the fact that $\hat {\cal A}$ is a finite set.

\noindent
We now need a Lemma to proceed with the proof. It is essentially a remark about the knot insertion formula in the one-dimensional case.

\smallskip
\begin{lemma}\label{C12sep1}
Let ${\bf \Sigma}=\{s_1,...,s_{n+p}\}$ be a knot vector, with the corresponding functions vectors ${\bf \Omega}_u=\{u_1,...,u_{n+p-1}\}$ and ${\bf \Omega}_v=\{v_1,...,v_{n+p-1}\}$. If we denote by $N_{h}^{(p)}(s)$ and $\bar N_{k}^{(p)}(s)$ the GB-splines of order $p$ constructed, respectively, on ${\bf \Sigma}$, ${\bf \Omega}_u$, ${\bf \Omega}_v$ and on the corresponding vectors obtained by inserting in ${\bf \Sigma}$, ${\bf \Omega}_u$, ${\bf \Omega}_v$ a certain number of knots and functions, then for any couple of different indices $1\le i,j\le n$, we get
\begin{align*}
&N_{i}^{(p)}(s)=\sum_{\bar \i_1\le h\le\bar \i_2}a_{h,p}\bar N_{h}^{(p)}(s), \notag \\
&N_{j}^{(p)}(s)=\sum_{\bar \j_1\le k\le\bar \j_2}b_{k,p}\bar N_{k}^{(p)}(s), %\label{17set3}
\end{align*}
where the $a_{h,p}$ and $b_{k,p}$ are coefficients obtained by a repeated application of the knot insertion formula. Moreover, we have either $\bar \i_2< \bar \j_2$ or $\bar \j_2<\bar \i_2$.
\end{lemma}

\smallskip
\noindent
{\bf Proof.} The Lemma is a direct consequence of the knot insertion formula. \hfill$\square$

\smallskip
If we assume that the column corresponding to $\bar A$ has two non-zero entries, then the two anchors $A_1,A_2\in {\cal A}$ corresponding to them, by Lemmas \ref{C17set1} and \ref{C12sep1}, are not right-horizontally and not down-vertically shifted, which contradicts the assumptions of the Theorem. In fact, if they are either right-horizontally or down-vertically shifted, one of the following possibilities occurs and for each of them it's not difficult to show, by using Lemma \ref{C17set1}, Proposition\ref{1mar14} and Lemma \ref{C12sep1}, that we can construct $\hat A\in\hat {\cal A}$ satisfying neither \eqref{129a} nor \eqref{129c}:

\begin{itemize}
\item $\bar A_1$ and $\bar A_2$ are right-horizontally shifted: it's possible to construct $\hat A\in\hat {\cal A}$ such that
\begin{equation*}
i^s_{p+1}(\hat A)> i^s_{p+1}(\bar A)\, \qquad i^t_{q+1}(\hat A)= i^t_{q+1}(\bar A),
\end{equation*}
which contradicts the definition of $\bar A$ given in \eqref{129a}-\eqref{129c}. Then, $\bar A_1$ and $\bar A_2$ cannot be right-horizontally shifted. 
\item $\bar A_1$ and $\bar A_2$ are down-vertically shifted: it's possible to construct $\hat A\in\hat {\cal A}$ such that
\begin{equation*}
i^t_{q+1}(\hat A)< i^t_{q+1}(\bar A),
\end{equation*}
which contradicts the definition of $\bar A$ given in \eqref{129a}-\eqref{129c}. Then, $\bar A_1$ and $\bar A_2$ cannot be down-vertically shifted. 
\end{itemize}

\noindent
Thus, we conclude that the column corresponding to $\bar A$ has only one non-zero entry.\hfill $\square$

\smallskip
\noindent
Finally, we show that analysis-suitable/dual-compatible T-meshes are always weakly dual-compatible.
\medskip
\begin{theorem}\label{NN12mag4}
A dual-compatible T-mesh $M$ is weakly dual-compatible (of any type).
\end{theorem}

\smallskip
\noindent
{\bf Proof.} Let $A_1,A_2\in {\cal A}\subset {\cal A}_{p,q}(M)$ and, without loss of generality, assume that they partially overlap horizontally with ${\bf I_l^s}(A_1)\cap{\bf I_l^s} (A_2) \neq \emptyset$. If ${\bf I_l^s}(A_1)\neq {\bf I_l^s}(A_2)$, then by applying the definition of partial horizontal overlapping, it can be proved that $A_1,A_2$ are left-horizontally and right-horizontally shifted. Otherwise, if ${\bf I_l^s}(A_1)={\bf I_l^2}(A_2)$, then $A_1$ and $A_2$ must overlap vertically, with ${\bf I_l^t}(A_1)\neq {\bf I_l^t}(A_2)$; then by applying the definition of partial vertical overlapping, $A_1,A_2$ are down-vertically and up-vertically shifted.  \hfill $\square$

\medskip
\noindent
Note that, as a consequence of Definition \ref{NN12mag2}, we get a characterization of the class of weakly dual-compatible T-meshes, which requires only the local vectors and the multiplicities of the knots. In order to give an example of weakly dual-compatible T-mesh which is not dual-compatible, let us consider the T-mesh $M$ drawn in Figure 7 with $p=q=\mu=\nu=4$ and all the knots in the active region having multiplicity $1$. It's easy to verify that the T-mesh of Figure 7 is weakly dual-compatible both of type RD and of type RU. By Definition \ref{C7} it's also evident that $M$ is not analysis-suitable. Therefore, we can state that the class of weakly dual-compatible T-meshes includes non-analysis-suitable elements.

\begin{figure}[h] 
\centering
\includegraphics[scale=0.44]{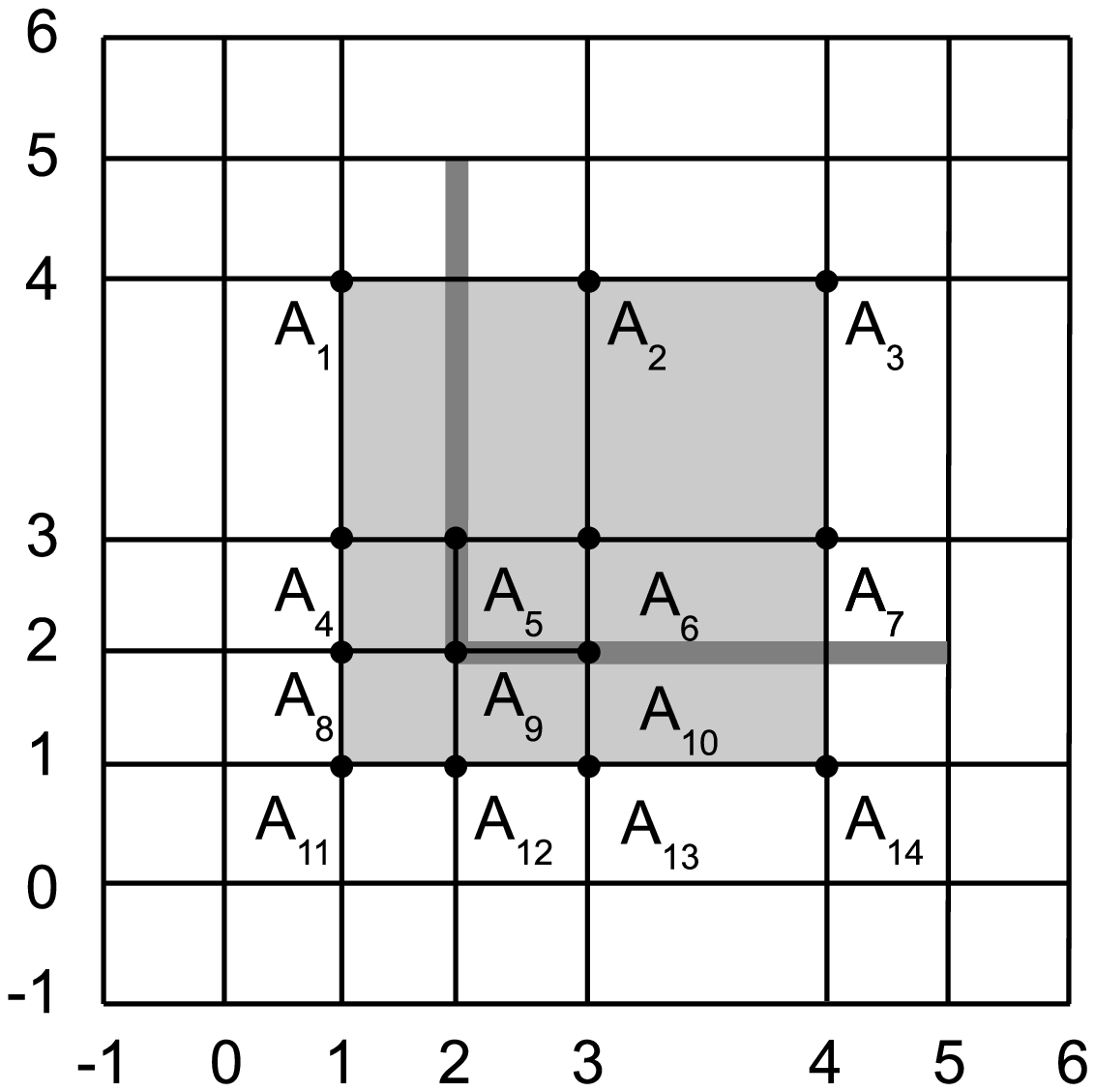}
\caption{Figure 7: {\sl A constructive example of weakly dual-compatible T-mesh and, which is clearly not analysis-suitable/dual-compatible (the T-junction extensions are highlighted in dark gray).}}
\end{figure}

\smallskip
\noindent
The characterization of weakly dual-compatible T-meshes is much simpler than that of VMCR T-meshes. This feature suggests the possibility to develop local refinement strategies preserving the weakly dual-compatible property. Let us give an example of refinement process based on weakly dual-compatible T-meshes, in order to show the potential use of this class of T-meshes.

\smallskip
\noindent
{\bf Example 4.} Let us consider the sequence of T-meshes in Figure 8 (we assume that the knots in the interior of the active region have multiplicity $1$). At each refinement step we add, alternatively, two horizontal or vertical edges subdividing the bottom/right-most $2\times 2$ square of cells of the active region. This strategy produces weak dual-compatible T-meshes: since the starting mesh is weakly dual-compatible of any type, always refining the bottom/right-most cells makes the anchors either left-horizontally or up-vertically shifted, and then the resulting T-mesh is always weakly compatible of type LU.
\begin{figure}[h] 
\centering
\includegraphics[scale=0.40]{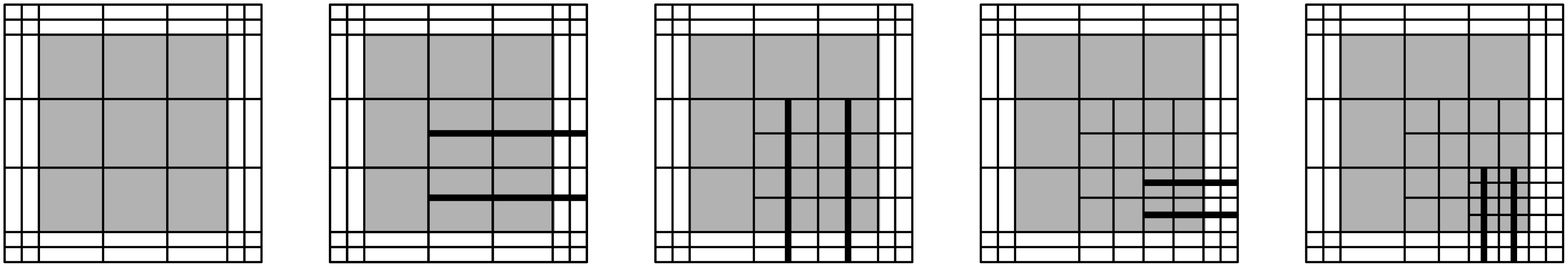}
\caption{Figure 8: {\sl An example of refinement producing weakly dual-compatible T-meshes, which are not dual-compatible. The active region is highlighted in grey, and bold lines represent the inserted edges.}}
\end{figure}

\section{Conclusions and future works}

\medskip
In this paper we completed the work of generalizing the T-spline approach to the GB-splines (started in \cite{bracco12b} for the trigonometric case). The obtained GT-splines are a flexible tool suitable to exactly model shapes otherwise impossible to exactly reproduce. The study of their linear independence led to the concept of VMCR T-meshes and of weakly dual-compatible T-meshes, two classes of T-meshes guaranteeing the linear independence of both the T-spline and GT-spline associated blending functions with the same bi-order. We proved, for any bi-order, that the classical dual-compatible/analysis-suitable T-meshes are weakly dual-compatible, and that, in turn, weakly dual-compatible T-meshes are VMCR. Moreover, we showed that there exist weakly dual-compatible T-meshes which are not analysis-suitable.

\smallskip
\noindent
GT-splines are different from classical T-splines: the local and global index vectors associated to the anchors determine not only the corresponding knot vectors but also the corresponding function vectors. The presence of the classes of VMCR T-meshes, which guarantees linear independence both for T-splines and for GT-splines, suggest an application of GT-spline to isogeometric analysis within the same framework used for T-splines. Note that the refinement strategies available for T-splines (see \cite{scot}) can be applied to GT-splines, since an analysis-suitable T-mesh is weakly dual-compatible. Moreover, the simple characterization of the new class of weakly dual-compatible T-meshes, suggests the possibility to develop new local refinement algorithms based on this larger class of T-meshes, which could be employed for the polynomial T-splines too.

\section*{Acknowledgements}
Durkbin Cho was supported by the Dongguk University Research Fund of 2012.

%%%%%%%%%%%%%%%%%%%%%%%%%%%%%%%%%%%%%%%%%%%%%%%%%%%%%%%%%

%%%%%%%%%%%%%%%%%%%%%%%%%%%%%%%%%%%%%%%%%%%%%%%%%%%%%%%%%

\end{document}